\documentclass{article}
\usepackage{graphicx} 

\usepackage{fullpage}

\usepackage{titling}

\usepackage{mathUtil}

\usepackage{fontawesome5}
\usepackage{xcolor}
\usepackage[outline]{contour}

\usepackage[backend=biber,style=numeric,maxnames=50]{biblatex}
\usepackage{hyperref}
\hypersetup{
    colorlinks=true,
    linkcolor=blue,
    filecolor=magenta,
    urlcolor=cyan,
    pdftitle={Colimits and Free Constructions of Heyting Algebras through Esakia Duality},
    pdfpagemode=FullScreen,
}
\title{Colimits and Free Constructions of Heyting Algebras through Esakia Duality}
\author{Rodrigo Nicolau Almeida}
\date{\today}

\begin{document}

\maketitle

\begin{abstract}
    In this paper\footnote{This work was originally developed in late 2024. It has suffered numerous changes over the last couple of years, having been expanded in several directions, with several results being added, corrected and reformatted. The present version serves to fix some minor typos and incorrections of the previous uploaded version. A fully reworked presentation of the results from Sections \ref{Free Constructions of Heyting Algebras over Distributive Lattices} and \ref{Applications} was given in \cite[Chapter 3]{almeidaphdthesis}, whilst the remaining sections were presented in a reformulated fashion in \cite[Chapter 4]{almeidaphdthesis}. These will be submitted as publications in their own right, as they contain substantial new material not present in this version. In both cases, the present preprint contains a lot of proofs of material of peripherical interest to the central constructions outlined in those chapters. In order to keep all such works consistent as much as possible, we have thus refrained from editing this beyond the strictly necessary.} we provide an explicit construction of the left adjoint to the inclusion of Heyting algebras in the category of distributive lattices using duality, generalizing a technique due to Ghilardi. This is employed to give a new concrete description of colimits and free Heyting algebras. Using such tools, we obtain in a uniform way some results about the category of Heyting algebra. These are both known, concerning the amalgamation property and related facts; and new, such as the fact that the category of Heyting algebras is comonadic over the category of distributive lattices. We also study some generalizations and variations of this construction to different related settings: (1) we analyse some subvarieties of Heyting algebras -- such as Boolean algebras, $\mathsf{KC}$ and $\mathsf{LC}$ algebras, and show how the construction can be adapted to such cases, and deriving some logical properties of $\mathsf{LC}$ in the process; (2) we study an adjunction between Heyting algebras and Boolean algebras, showing how this provides a categorical semantics for inquisitive logic; (3) we study the relationship between the category of image-finite posets with p-morphisms and the category of posets with monotone maps, and using a variation of the above ideas, provide an adjunction between such categories, with applications in the coalgebraic approach to intuitionistic modal logic.
\end{abstract}

\tableofcontents

\section{Introduction}

Free algebras -- in their form of \textit{Lindenbaum-Tarski algebras} -- play an important role in analyzing systems of non-classical logics \cite{Chagrov1997-cr,Blackburn2002-fd,Galatos2007-vu}. For this purpose, \textit{duality theory} has played a remarkable role, turning the abstract ``word constructions" of free algebras into concrete and analyzable objects: for instance, to calculate coproducts of Boolean algebras, one need only describe products of Stone spaces, which turn out to be exactly the usual topological products. Since arbitary free algebras can be computed as coproducts of $1$-generated algebras, this also allows an easy computation of free Boolean algebras. A similar and more general situation holds concerning distributive lattices and their prime spectra.

In this paper we focus on the case of \textit{Heyting algebras}, which are the algebraic models of $\mathsf{IPC}$, the intuitionistic propositional calculus, or intuitionistic logic. Despite having a good representation theory in terms of Esakia spaces, the structure of free Heyting algebras is much more difficult than that of Boolean algebras or distributive lattices: already the free algebra on one generator is infinite, observed by Rieger \cite{Rieger1952-RIEOTL} and Nishimura \cite{nishimuralattice}. The free algebra on two generators is notoriously difficult, having no, known, easy lattice theoretic description. 

As a result, considerable attention has been devoted to describing such objects through recourse to more combinatorial and topological structures. Key results in this direction include:
\begin{enumerate}
    \item The construction of the $n$-universal model (due independently to Bellissima \cite{Bellissima1986}, Grigolia \cite{Grigolia1995} and Shehtman \cite{shehtmanrieger}, see \cite[8.8]{Chagrov1997-cr} or \cite[3.2]{bezhanishviliphdthesis} for modern presentations
    ), which clarified the size and structure of irreducible elements of the free $n$-generated algebra, as well as some of its properties like completeness;
    \item The construction of the free $n$-generated Heyting algebra by Ghilardi \cite{ghilardifreeheyting} (which generalized previous work by Urquhart \cite{Urquhart1973}, and was in turn generalized for all finitely presented Heyting algebras by Butz \cite{Butz1998}), which allowed to show that such algebras in fact carry more structure, like a dual connective to the implication.
\end{enumerate}

Despite the importance of such advances, the situation presently in the literature presents a big contrast between what is known about Boolean algebras, distributive lattices and Heyting algebras: the former have descriptions of their free algebras, coproducts, and free one-step extensions for any number of generators, and any original algebra, whilst the latter are only available for special cases, namely those involving finite algebras.

The main contribution of this paper is to provide a generalization of the constructions presented by Ghilardi to the infinite case. The main hurdle in this is realizing wherein finiteness plays a role, to allow using Priestley duality in place of Birkhoff duality for finite distributive lattices, a task which turns out to not be obvious.  This is equivalent to describing, in dual terms, the left adjoint to the inclusion of $\mathbf{HA}$, the category of Heyting algebras, in $\mathbf{DL}$, the category of distributive lattices. Due to the duality theory for such categories, this yields a right adjoint to the inclusion of $\mathbf{Esa}$, the category of Esakia spaces and continuous p-morphisms, into the category $\mathbf{Pries}$ of Priestley spaces and continuous order-preserving maps. Having this at hand, we employ it for several purposes:
\begin{enumerate}
    \item We provide a description of the free Heyting algebras on any number of generators, and provide a concrete description of coproducts of Heyting algebras as well as of the pushout of Heyting algebras. From this we derive some results concerning co-distributivity of Heyting algebras.
    \item We describe explicitly the adjunction between Boolean algebras and Heyting algebras concerning the regularization functor $\neg\neg$. This provides categorical semantics for inquisitive logic, connecting to the algebraic and topological semantics previously given to it.
    \item We first show that this construction can be adapted when dealing extensions of $\mathsf{IPC}$, or algebraically, subvarieties of $\mathbf{HA}$. Rather than providing a general theory, we focus on three important and illustrative cases: Boolean algebras, $\mathbf{KC}$ algebras and G\"{o}del algebras. As a consequence of our analysis we show that $\mathsf{LC}$ has \textit{uniform local tabularity} (first studied in \cite{Shehtman2016}under the name \textit{finite formula depth property}), and that in fact every formula is equivalent to one of implication depth $2$.
    \item Finally, we look at the relationship between the category of posets with p-morphisms and the category of posets with monotone maps, and provide a right adjoint to the inclusion which is heavily inspired by the above constructions. The latter has found applications in the coalgebraic study of intuitionistic modal logic \cite{coalgebraicintuitionisticmodal}, addressing the open problem of how to represent Kripke frames for intuitionistic modal logics coalgebraically. We also show that the two constructions can be connected by the concept of an order-compactification.
\end{enumerate}

The structure of the paper is as follows: in Section \ref{Preliminaries} we recall some basic preliminaries, and fix some notation moving forward. Section \ref{Free Constructions of Heyting Algebras over Distributive Lattices} contains the main technical tools, in the form of the description of the free Heyting algebra generated by a distributive lattice preserving a given set of relative pseudocomplements; in Section \ref{Applications} we showcase some basic applications of the aforementioned theory to the study of Heyting algebras: some well-known results, like the amalgamation theorem, are shown here through direct methods, as well as some seemingly new results, like the comonadicity of the category of Heyting algebras over the category of distributive lattices. In Section \ref{Regular Heyting algebras and free Heyting extensions of Boolean algebras} we look at an  adjunction holding between Boolean algebras and Heyting algebras, related to the regularization functor, and show that the left adjoint of such a functor admits a description similar to the one given by the dual construction just given, and from this some applications to the study of inquisitive logic. In Section \ref{Subvarieties of Heyting algebras} we look at the case of subvarieties of Heyting algebras and illustrate how these constructions can be fruitfully analysed there. In Section \ref{The Category of Posets with P-morphisms} we look at the corresponding relationship between the category of posets with p-morphisms and the category of posets with monotone maps. We conclude in Section \ref{Conclusion} with some outlines of future work.

\section{Preliminaries}\label{Preliminaries}

We assume the reader is familiar with the basic theory of Boolean algebras and Distributive lattices, as well as with Stone duality (see for example \cite{Davey2002-lr}). Throughout we will use the term ``distributive lattices" for \textit{bounded} distributive lattices, i.e., our lattices will always have a bottom element $(0)$ and a top element $(1)$. We also assume the reader is familiar with Heyting algebras and their elementary theory (see e.g. \cite{Esakiach2019HeyAlg}).

Key to our results will be the use of duality theory for distributive lattices and Heyting algebras, in the form of Priestley and Esakia duality. Recall that an \textit{ordered topological space} is a triple $(X,\leq,\tau)$ where $(X,\leq)$ is a poset and $(X,\tau)$ is a topological space. Often we refer to $(X,\leq)$, or just $X$ as an ordered topological space, when no confusion arises. Given such a space, and a subset $U$, we write ${\uparrow}U=\{x : \exists y\in U, x\geq y\}$ (symmetrically, ${\downarrow}U$), and say that $U$ is an \textit{upset} (downset) if $U={\uparrow}U$ ($U={\downarrow}U$).

\begin{definition}
    Let $X$ be an ordered topological space. We say that $X$ is a \textit{Priestley space} if $X$ is compact and whenever $x\nleq y$ there is some clopen upset $U$ such that $x\in U$ and $y\notin U$.

    We say that a Priestley space $X$ is \textit{Esakia} if whenever $U$ is clopen, ${\downarrow}U$ is clopen as well.
\end{definition}

\begin{definition}
    Let $X,Y$ be Priestley spaces. We say that $f:X\to Y$ is a \textit{p-morphism} if whenever $f(x)\leq y$ there is some $x'\geq x$ such that $f(x')=y$.
\end{definition}

We denote by $\mathbf{DLat}$ the category of distributive lattices and lattice homomorphisms, and by $\mathbf{Pries}$ the category of Priestley spaces and monotone and continuous functions. We denote by $\mathbf{HA}$ the category of Heyting algebras and Heyting algebra homomorphisms, and by $\mathbf{Esa}$ the category of Esakia spaces and p-morphisms.

\begin{theorem}
    The category $\mathbf{DLat}$ is dually equivalent to the category $\mathbf{Pries}$, and the category $\mathbf{HA}$ is dually equivalent to the category $\mathbf{Esa}$.
\end{theorem}

The constructions in this duality are similar to those of Stone duality: one takes for a given distributive lattice or Heyting algebra $\alg{H}$ the spectrum $\mathsf{Spec}(H)=\{x : x \text{ is a prime filter in $H$}\}$; and conversely, given a Priestley or Esakia space $X$, one takes $\mathsf{ClopUp}(X)$, the set of clopen upsets. More specifically, given a distributive lattice $\mathcal{D}$, we denote the representation map as
\begin{equation*}
    \phi:\alg{D}\to \mathsf{ClopUp}(X_{D})
\end{equation*}
where $\phi(a)=\{x\in X_{D} : a\in x\}$. The reason this map is bijective follows from the prime ideal theorem: given $\alg{D}$ a distributive lattice, and $A\subseteq D$, we write $\mathsf{Fil}(A)$ for the filter generated by $A$, and similarly, $\mathsf{Id}(A)$ for the ideal generated by $A$. We also recall that the following prime filter theorem holds for distributive lattices:

\begin{theorem}\label{Prime Ideal Theorem for Distributive Lattices}
    Let $\alg{D}$ be a distributive lattice, $F,I\subseteq D$ be respectively a filter and an ideal, such that $F\cap I=\emptyset$. Then there exists a prime filter $F'\supseteq F$ such that $F'\cap I=\emptyset$.
\end{theorem}

As far as morphisms work, given a distributive lattice homomorphism $f:D\to D'$, the dual morphism is given by $f^{-1}:X_{D'}\to X_{D}$, which maps prime filters to prime filters. Similarly, if $p:X\to Y$ is a monotone continuous map, then $p^{-1}:\mathsf{ClopUp}(Y)\to \mathsf{ClopUp}(X)$ is a distributive lattice homomorphism; this restricts to Esakia spaces and Heyting algebras in the obvious way.

A fact that we will often use below is that if $(X,\leq)$ is a Priestley space, and $Y\subseteq X$ is a closed subset, then $(Y,\leq_{\restriction Y})$ is a Priestley space as well. We also note the following important, though easy to prove, fact:

\begin{proposition}\label{Specific Relative pseudocomplements}
    Let $\alg{D}$ be a distributive lattice, and 
let $a,b\in \alg{D}$ be arbitrary. Then if $c$ is the relative pseudocomplement of $a$ and $b$, then
\begin{equation*}
    \phi(c)=X_{D}-{\downarrow}(\phi(a)-\phi(b)).
\end{equation*}
\end{proposition}

Given an arbitrary subset $U\subseteq X$ of an ordered topological space, we will then write $\Box U=X-{\downarrow}X-U$.

We also recall the following construction, which is part of the folklore of the subject (see \cite{gehrkevangoolbooktopologicalduality} for an in-depth discussion):

\begin{definition}\label{Description of free algebras}
    Let $X$ be an arbitrary set. Let $\mathbf{D}(X)$ be the distributive lattice dual to the Priestley space
    \begin{equation*}
        F_{D}(X)\coloneqq \mathbf{2}^{X}
    \end{equation*}
    where $\mathbf{2}$ is poset $0<1$, and $F_{D}(X)$ is given the product topology and the product ordering.
\end{definition}

We then have the following well-known fact:

\begin{proposition}
    If $X$ is an arbitrary set, then $\mathbf{D}(X)$ is isomorphic to the free distributive lattice on $X$ many generators. Moreover, if $X$ is finite, then $F_{D}(X)$ is isomorphic to  $(\mathcal{P}(X),\supseteq)$with the discrete topology.
\end{proposition}

We will also occasionally need the \textit{Rieger-Nishimura ladder}, which is the Esakia dual of the $1$-generated Heyting algebra. This can be constructed (see e.g. \cite{Bezhanishvili2006}) in the following way: starting with three points $2,1,0$, and order them by setting $2<1$, and let $X_{0}=\{2,1,0\}$. We then construct $X_{n+1}$ by picking for each antichain of points in $X_{n}$ $x_{0},...,x_{n}$ such that there is no point $y$ immediately below each $x_{i}$; as it turns out, one can show that $X_{n+1}$ always adds only two points. We let $X(1)=\bigcup_{n\in \omega}X_{n}\cup {\infty}$, and order $\infty$ below every other point; we then give $X(1)$ the topology based on finite sets and cofinite sets containing $\infty$. The resulting space is depicted in Figure

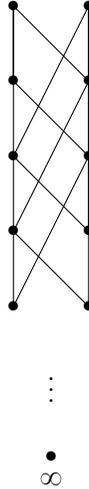
\begin{figure}[h]
    \centering
\begin{tikzpicture}
    \node at (0,0) {$\bullet$};
    \node at (1,0) {$\bullet$};
    \node at (0,-1) {$\bullet$};
    \node at (1,-1) {$\bullet$};
    \node at (0,-2) {$\bullet$};
    \node at (1,-2) {$\bullet$};
    \node at (0,-3) {$\bullet$};
    \node at (1,-3) {$\bullet$};
    \node at (0,-4) {$\bullet$};
    \node at (1,-4) {$\bullet$};
    \node at (0.5,-5) {$\vdots$};
    \node at (0.5,-6) {$\bullet$};
    \node at (0.5,-6.3) {$\infty$};

    \draw (0,0) -- (0,-1) -- (0,0) -- (1,-1) -- (1,0) -- (1,-1) -- (1,-2) -- (0,-1) -- (0,-2) -- (1,0);

    \draw (0,-2) -- (0,-3) -- (1,-1) -- (1,-2) -- (1,-3) -- (0,-2);

    \draw (0,-3) -- (0,-4) -- (1,-2) -- (1,-3) -- (1,-4) -- (0,-3);
\end{tikzpicture}
    
    \caption{Rieger-Nishimura Ladder}
    \label{fig:riegernishimura}
\end{figure}

An important part of our work will deal with inverse limits of Priestley spaces. The following is clear from duality theory:

\begin{proposition}\label{Directed Unions are dual to inverse limits}
    Let $(\alg{D}_{n},\iota_{i,j})$ be a chain of distributive lattices, connected by embeddings $i_{i,j}:\alg{D}_{i}\to \alg{D}_{j}$ for $i<j$ satisfying the usual compatibility laws of directed systems. Let $(X_{n},p_{i,j})$ be the inverse system of Priestley spaces one obtains by dualising all arrows. Then the directed union of the directed chain, $\alg{D}_{\omega}$, is dual to the projective limit $X_{\omega}$ of the inverse system. Moreover, if each $X_{n}$ is an Esakia space, and each $p_{i,j}$ is a p-morphism, then $X_{\omega}$ is likewise an Esakia space.
\end{proposition}
\begin{proof}
    Note that if each $X_{n}$ is an Esakia space, and each $p_{i,j}$ is a p-morphism, then the inverse limit is dual to the directed limit of distributive lattices; in these conditions, one can define a Heyting algebra structure on the directed union by taking the equivalence classes of the relative pseudocomplements arising in the finitely many stages, which will be well-defined since $p_{i,j}^{-1}$ are Heyting algebra homomorphisms.
\end{proof}

We also recall the well-known construction of the \textit{Vietoris hyperspace}, due to Leopold Vietoris \cite{Vietoris1921} (see also \cite{guramremarksonhyperspacesforpriestleyspaces}). This can equivalently be given as follows: if $(X,\leq)$ is a Priestley space, let $V(X)$ be the set of all non-empty closed subsets of $X$. We give this set a topology by considering a subbasis consisting of, for $U,V\subseteq X$ clopen sets:
\begin{equation*}
    [U]=\{C\in V(X) : C\subseteq U\} \text{ and } \langle V\rangle =\{C\in V(X) : C\cap V\neq\emptyset\}.
\end{equation*}
This is sometimes called the ``hit-and-miss topology", and the resulting space is called the \textit{Vietoris hyperspace}\footnote{Vietoris' original definition considers instead $[U]$ and $\langle V\rangle$ for $U$ open. But note that these topologies are equivalent when $X$ is Stone, since closed subsets are compact in such spaces.}. On this space we define an order relation $C\preceq D$ if and ony if\footnote{The reason this order is reverse inclusion has to do with the choice to use \textit{upsets} rather than downsets; in \cite{ghilardifreeheyting}, the author uses downsets, which is why the swap no reversion is needed.} $D\subseteq C$. Then we have the following:

\begin{lemma}\label{Vietoris is Priestley}
    The space $(V(X),\preceq)$ is a Priestley space. Indeed, it is an Esakia space.
\end{lemma}
\begin{proof}
    The fact that $V(X)$ is compact is a standard fact (see e.g. \cite{Engelking1989-qo}), but we prove it here for completeness. Assume that $V(X)=\bigcup_{i\in I}[U_{i}]\cup \bigcup_{j\in J}\langle V_{j}\rangle$ is a cover by clopen sets. Now consider $C=X-\bigcup_{i\in I}V_{i}$. If $C$ is empty, then the sets $V_{i}$ cover $X$ and so finitely many of them cover $X$, say $V_{i_{0}},...,V_{i_{n}}$. Then if $A\in V(X)$, then $A$ is non-empty and so contains some $x$ which must lie in some $V_{i_{k}}$ for $k\in \{1,...,n\}$, and so $A\in \langle V_{i_{k}}\rangle$. So the sets $\{\langle V_{i_{k}}\rangle : k\in \{1,...,n\}\}$ form a finite subcover. 

Otherwise, we have that $C\neq \emptyset$, so $C\in V(X)$, and since $C\notin \langle V_{j}\rangle$ for each $i$, by definition $C\in [U_{i}]$ for some $i$. Since $X-\bigcup_{i\in I}V_{i}\subseteq U_{i}$ then $X-U_{i}\subseteq \bigcup_{i\in I}V_{i}$, so by compactness, $X-U_{i}\subseteq V_{i_{0}}\cup...\cup V_{i_{n}}$. Hence consider the finite subcover $\{[U_{i}]\}\cup \{\langle V_{i_{j}}\rangle : j\in \{1,...,n\}\}$. If $D$ is any closed set, and $D\subseteq U_{i}$, then we are done; otherwise, $D\cap X-U_{i}\neq \emptyset$, so $D\cap V_{i_{0}}\cup...\cup V_{i_{n}}\neq \emptyset$, and so the result follows.

Now we prove the Priestley Separation axiom. Assume that $C\npreceq D$ are closed subsets. Since the space $X$ is a Stone space, we know that $C=\bigcap_{i\in I}V_{i}$ and $D=\bigcap_{j\in J}W_{j}$ where $V_{i},W_{j}$ are clopen. Hence there must be some $i$ such that $C\subseteq V_{i}$ and $D\nsubseteq V_{i}$; i.e., $C\in [V_{i}]$ and $D\notin [V_{i}]$. Since $[V_{i}]$ is clearly a clopen upset under this order, we verify the Priestley separation axiom.

Finally, consider an arbitrary clopen difference:
    \begin{equation*}
        Z=[U]\cap \langle V_{0}\rangle\cap...\cap \langle V_{n}\rangle.
    \end{equation*}
    Then we want that ${\downarrow}Z$ is clopen. Indeed consider:
    \begin{equation*}
        \bigcap_{i=1}^{n}\langle U\cap V_{i}\rangle.
    \end{equation*}
    If $C$ belongs to this set, then it contains a subset $D=\{x_{1},...,x_{n}\}$, elements in $V_{i}$ and $U$; so $D\subseteq U$, it is a closed subset, and intersects each of the sets. So $C\in {\downarrow}Z$. Conversely, if $C\in {\downarrow}Z$, then it is clear it belongs above. So since $\langle U\cap V_{i}\rangle$ are all subbasic clopens, we get that ${\downarrow}Z$ is a clopen, as desired.
\end{proof}

\section{Free Constructions of Heyting Algebras over Distributive Lattices}\label{Free Constructions of Heyting Algebras over Distributive Lattices}

In this section we explain how one can construct a Heyting algebra freely from any given distributive lattice. This construction generalises the case of Heyting algebras freely generated by finite distributive lattices, as originally analyzed by Ghilardi \cite{ghilardifreeheyting}.

\subsection{Conceptual Idea of the Construction and the Finite Case}\label{Conceptual idea}

Before diving into the details, we will provide an informal explanation of the construction by Ghilardi which justifies the technical developments, following the core of the discussion from \cite{Bezhanishvili2011}. We urge the reader to consult this section for intuition whilst going through the details of the next section.

Suppose that $D$ is a finite distributive lattice, and that we wish to construct a Heyting algebra from $D$. Then certainly we must (1) freely add to $D$ implications of the form $a\rightarrow b$ for $a,b\in D$. This is the same as freely generating a distributive lattice out of $D^{2}$, and considering its coproduct with $D$. Dually, if $X$ is the dual poset to $D$, this will then be the same as considering
\begin{equation*}
    X\times \mathcal{P}(X\times X),
\end{equation*}
given that the free distributive lattice on a set of pairs of generators $X$ is dual to $\mathcal{P}(X\times X)$ (see Section \ref{Preliminaries}). However, we need (2) to impose axioms forcing these implications to act like relative pseudocomplements. A first move is to impose the axiom of a Weak Heyting Algebra:
\begin{enumerate}
    \item $a\rightarrow a=1$;
    \item $(a\vee b)\rightarrow c=(a\rightarrow c)\wedge (b\rightarrow c)$;
    \item $a\rightarrow (b\wedge c)=(a\rightarrow b)\wedge (a\rightarrow c)$;
    \item $(a\rightarrow b)\wedge (b\rightarrow c)\leq a\rightarrow c$.
\end{enumerate}
Applying a quotient under these axioms, will dually yield that from our $X$, we obtain (see \cite[Theorem 3.5]{Bezhanishvili2011}),
\begin{equation*}
    \mathcal{P}(X),
\end{equation*}
Now, reformulating this slightly, one can see that the upsets of $\mathcal{P}(X)$ are of the form
\begin{equation*}
    [U]=\{C\subseteq X : C\subseteq U\}
\end{equation*}
for $U$ a subset of $X$. This then provides an expansion of our lattice, since we can consider the map:
\begin{align*}
    i_{0}:D &\to \mathsf{Up}(\mathcal{P}(X))\\
    a &\mapsto [\phi(a)],
\end{align*}
which is easily seen to be injective, and a meet-homomorphism preserving the bounds.

The fact that we want to obtain a genuine relative pseudocomplement, means we need to impose further axioms, which implies throwing out some of the elements from $\mathcal{P}(X)$. The additional fact that we would want the map $i_{0}$ to be a distributive lattice homomorphism suggests a way of doing this: make it so that for each pair of upsets $U,V$, $[U]\cup [V]=[U\cup V]$. After some thought one is lead to consider \textit{rooted} subsets of $X$, obtaining the poset $\mathcal{P}_{r}(X)$. Some verifications show that this will indeed be the poset that is needed: that
\begin{equation*}
    i_{0}:D\to \mathsf{Up}(\mathcal{P}_{r}(X))
\end{equation*}
will be a distributive lattice embedding as desired, and that it will contain relative pseudocomplements for the elements from $D$, namely, for $a,b\in D$, the element $[-\phi(a)\cup \phi(b)]$.

Now at this point we will have added all implications to elements $a,b\in D$, obtaining a distributive lattice $D_{1}$, but all the new implications added might not in turn have implications between themselves. So we need to (3) iterate the construction, infinitely often, to add all necessary implications. However, the final complication is that each step of this construction adds implications to every element in the previous lattice. So in particular, $D_{2}$ will contain a fresh relative pseudocomplement, $[a]\rightarrow [b]$, which need not agree with the relative pseudocomplement $[a\rightarrow b]$. If we let the construction run infinitely often in this way, it could then be that in the end no element would be the relative pseudocomplement of $a$ and $b$, so we need to ensure that on the second iteration, the previously added relative pseudocomplements are preserved, i.e.:
\begin{equation*}
    i_{1}(i_{0}(a)\rightarrow i_{0}(b))=i_{1}(i_{0}(a))\rightarrow i_{1}(i_{0}(b))
\end{equation*}

In other words, we need the map $i_{1}$ to preserve the relative pseudocomplements of the form $i_{0}(a)\rightarrow i_{0}(b)$. It is this need which justifies the notion of a $g$-open subset detailed below, and leads to us considering, at last, as our one step construction, the poset
\begin{equation*}
    \mathcal{P}_{g}(X)
\end{equation*}
where $g$ is some order-preserving map which serves to index the relative pseudocomplements which are to be preserved.

To conclude this section, let us take a look at an example calculation:

\begin{example}\label{step by step construction of Rieger-Nishimura}
    Let $X=\mathbf{2}$ the $0<1$ poset. Then we can illustrate the first four steps of the step-by-step construction of Rieger-Nishimura ladder (see Figure \ref{fig:riegernishimura}) in Figure \ref{fig:stepstoconstructRiegerNishimura}. Note that the step-by-step construction does not in general agree with the construction we provided before.

    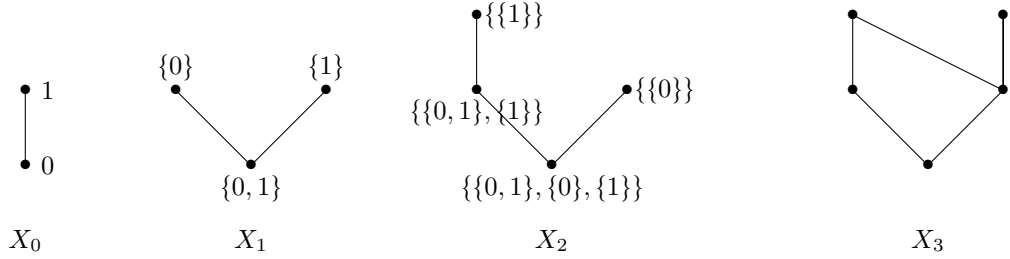
\begin{figure}[h]
        \centering
\begin{tikzpicture}
    \node at (0,0) {$\bullet$};
    \node at (0.3,0) {$0$};
    \node at (0,1) {$\bullet$};
    \node at (0.3,1) {$1$};
    \node at (0,-1) {$X_{0}$};

    \draw (0,0) -- (0,1);

    \node at (2,1) {$\bullet$};
    \node at (2,1.3) {$\{0\}$};
    \node at (3,0) {$\bullet$};
    \node at (3,-0.3) {$\{0,1\}$};
    \node at (4,1) {$\bullet$};
    \node at (4,1.3) {$\{1\}$};
    \node at (3,-1) {$X_{1}$};

    \draw (2,1) -- (3,0) -- (4,1);

\node at (6,2) {$\bullet$};
\node at (6.5,2) {$\{\{1\}\}$};
    \node at (6,1) {$\bullet$};
    \node at (6,0.7) {$\{\{0,1\},\{1\}\}$};
    \node at (7,0) {$\bullet$};

    \node at (7,-1) {$X_{2}$};

    \node at (7,-0.3) {$\{\{0,1\},\{0\},\{1\}\}$};
    \node at (8,1) {$\bullet$};
    \node at (8.5,1) {$\{\{0\}\}$};

    \draw (6,2) -- (6,1) -- (7,0) -- (8,1);

    \node at (11,1) {$\bullet$};
    \node at (11,2) {$\bullet$};
    \node at (13,2) {$\bullet$};
    \node at (13,1) {$\bullet$};
    \node at (12,0) {$\bullet$};
    \node at (12,-1) {$X_{3}$};

    \draw (11,1) -- (11,2) -- (13,1) -- (13,2) -- (13,1) -- (12,0) -- (11,1);

\end{tikzpicture}        \caption{Four steps of the construction of the Rieger-Nishimura}
        \label{fig:stepstoconstructRiegerNishimura}
    \end{figure}
\end{example}

\subsection{G-Open Subsets and Vietoris Functors}

We begin by adapting the notion of $g$-openness to our setting. Let $X,Y$ be two Priestley spaces, and let $f:X\to Y$ be a monotone and continuous map. By Priestley duality, such maps are dual to distributive lattice homomorphisms $f^{-1}:\mathsf{ClopUp}(Y)\to \mathsf{ClopUp}(X)$.

\begin{definition}
Let $X,Y,Z$ be Priestley spaces, and $f:X\to Y$ and $g:Y\to Z$ be monotone and continuous functions. We say that $f$ is \textit{open relative to $g$} ($g$-open for short) if $f^{-1}$ preserves relative pseudocomplements of the kind $g^{-1}[U]\rightarrow g^{-1}[V]$ where $U,V\in \mathsf{ClopUp}(Z)$.
\end{definition}

The following is a condition equivalent to being $g$-open, expressed purely in terms of the order:

\begin{equation}
\tag{*} \forall a\in X, \forall b\in Y, (f(a)\leq b \implies \exists a'\in X, (a\leq a' \ \& \ g(f(a'))=g(b)).
\end{equation}

\begin{lemma}\label{G-open equivalent to preservation of pseudocomplements}
Given $f:X\to Y$ and $g:Y\to Z$, we have that $f$ is $g$-open if and only if $f$ satisfies condition (*).
\end{lemma}
\begin{proof}
Asssume that $f$ satisfies (*). We want to show that $f^{-1}[g^{-1}[U]\rightarrow g^{-1}[V]]=f^{-1}g^{-1}[U]\rightarrow f^{-1}g^{-1}[V]$. Note that the left to right inclusion always holds because $f$ is a monotone map. So assume that $x\in f^{-1}g^{-1}[U]\rightarrow f^{-1}g^{-1}[V]$. Suppose that $f(x)\leq y$, and $y\in g^{-1}[U]$. By assumption, there is some $x'$ such that $x\leq x'$ and $g(f(x'))=g(y)$; hence $f(x')\in g^{-1}[U]$, so by assumption, $x'\in f^{-1}g^{-1}[U]$, and so, $x'\in f^{-1}g^{-1}[V]$. This means that $f(x')\in g^{-1}[V]$, so $y\in g^{-1}[V]$. This shows what we want.

Now assume that $p=f^{-1}$ is $g$-open, and $q=g^{-1}$. 
Assume that $f(x)\leq y$ where $x\in X$ and $y\in Y$. By duality, and abusing notation, this means that $p^{-1}[x]\subseteq y$. Consider $\mathsf{Fil}(x,\{p(q(a)) : q(a)\in y\})$ and $\mathsf{Id}(\{p(q(b)) : q(b)\notin y\})$. We claim these do not intersect. Because otherwise, for some $c\in x$, $c\wedge p(q(a))\leq p(q(b))$. Hence $c\leq p(q(a))\rightarrow p(q(b))$, since these exist and are preserved by $p$, and so $c\leq p(q(a)\rightarrow q(b))$. So $p(q(a)\rightarrow q(b))\in x$. Hence $q(a)\rightarrow q(b)\in p^{-1}[x]$, and so $q(a)\rightarrow q(b)\in y$, a contradiction. Hence by the Prime filter theorem (see Theorem \ref{Prime Ideal Theorem for Distributive Lattices}), we can extend $x$ to a filter $x'$ which does not intersect the presented ideal. By definition, working up to natural isomorphism, we then have that $g(f(x'))=g(y)$, which was to show.
\end{proof}

Now given $g:X\to Y$ a monotone and continuous function, and $S\subseteq X$ a closed subset, we say that $S$ is \textit{g-open} (understood as a Priestley space with the induced order and topology) if the inclusion is itself $g$-open. This means by the above lemma that $S$ is $g$-open if the following condition holds:

\begin{equation*}
\forall s\in S, \forall b\in X (s\leq b \implies \exists s'\in S (s\leq s' \ \& \ g(s')=g(b)).
\end{equation*}

Following Ghilardi, this can be thought of as follows: if we think of $X$ as represented by fibers coming from $g$, then whenever ${\uparrow}s$ meets an element of any fiber, then ${\uparrow}s\cap S$ must actually contain an element of that fiber. With this intuition, it is not difficult to show that if $x$ is arbitrary and $S$ is $g$-open, then $S\cap {\uparrow}x$ is $g$-open as well.

\begin{remark}\label{Remark on multiple maps}
   The fact that we pick a given unique $g$ is totally incidental; in fact, whilst the results of this section will be proved for a single $g$, one could take any number of continuous morphisms and obtain the same results.  We will need this briefly in the description of coproducts made below.
\end{remark}

Throughout, fix $g:X\to Y$ a continuous and order-preserving map (the case for several such maps being preserved is entirely similar). Recall from Lemma \ref{Vietoris is Priestley} that $(V(X),\preceq)$ is a Priestley space. From this space we can move closer to the space we will be interested, by first considering $V_{r}(X)\subseteq V(X)$, the space of \textit{rooted} closed subsets, with the induced order and the subspace topology.
On this we can prove the following\footnote{This is stated without proof in \cite[Lemma 6.1]{Bezhanishvili2011heytingcoalgebra}. The key idea of the proof below was communicated to me by Mamuka Jibladze.}:

\begin{lemma}
The space $(V_{r}(X),\preceq)$ is a Priestley space.
\end{lemma}
\begin{proof}
It suffices to show that $V_{r}(X)$ is closed. For that purpose, consider the following two subsets of $V(X)\times X$:
\begin{equation*}
    LB(X)=\{(C,r) : C\subseteq {\uparrow}r\} \text{ and } \ni_{X}=\{(C,r) : r\in C\}.
\end{equation*}
We note that $\ni_{X}$ is closed\footnote{Indeed, one could simply observe that this is the pullback of the continuous map $\{\}:X\to V(X)$ along the second projection of $\leq\subseteq V(X)\times V(X)$.}: if $(C,x)$ is such that $x\notin C$, then $\{x\}\cap C=\emptyset$. Since $X$ is a Stone space, let $U$ be a clopen subset separating them i.e., $x\in U$ and $C\subseteq X-U$. Then we can consider,
\begin{equation*}
    [X-U]\times U
\end{equation*}
which is a neighbourhood of $(C,x)$; and certainly if $(D,y)\in [X-U]\times U$, then $y\notin D$.

In turn, $LB(X)$ can also be seen to be closed, through a direct argument: suppose that $(C,r)\notin LB(X)$. By definition, then, $C\nsubseteq {\uparrow}r$, so there is a point $y\in C$ such that $r\nleq y$. By the Priestley separation axiom, there is a clopen upset $U$ such that $r\in U$ and $y\notin U$. So consider
\begin{equation*}
    S=\langle X-U\rangle \times U.
\end{equation*}
This is clearly an open subset of $V(X)\times X$, and $(C,r)$ belongs there by hypothesis. But also, if $(D,k)\in S$, then $D-U\neq \emptyset$, so there is some point $m\notin U$; but since $U$ is an upset, then we must have $k\nleq m$, i.e., $(D,k)\notin LB(X)$. This shows that $LB(X)$ is closed.

Now note that $$V_{r}(X)=\pi_{V(X)}[LB(X)\cap \ni_{V(X)}].$$

Since $\pi$ is a closed map (given it is a continuous surjection between Stone spaces), we then have that $V_{r}(X)$ is a closed subspace of $V(X)$, as desired.\end{proof}

Finally, we will refine the rooted subsets to the ones we are truly interested in -- the $g$-open ones. Let $V_{g}(X)$ be the set of $g$-open, closd and rooted subsets.  The following lemma encapsulates the key fact about this.

\begin{lemma}\label{The g-open subspace is also a Priestley space}
Assume that $g:X\to Y$ is an order-preserving and continuous map, and that $X$ has relative pseudocomplements of sets of the form $g^{-1}[U]$. Then the subspace $V_{g}(X)\subseteq V_{r}(X)$ is a closed subspace, and hence is a Priestley space as well with the induced order.
\end{lemma}
\begin{proof}
Assume that $M$ is not $g$-open. This means that there is some $x\in M$ such that $x\leq y$ and for each $k\in M\cap {\uparrow}x$, $g(k)\neq g(y)$. Thus $g(y)\notin g[M\cap {\uparrow}x]$, so there is a clopen $-U\cup V$ such that $g[M\cap {\uparrow}x]\subseteq -U\cup V$ and $g(y)\in U-V$. Then $M\cap {\uparrow}x\subseteq g^{-1}[-U\cup V]$, and so $M\subseteq -{\uparrow}x \cup g^{-1}[-U\cup V]$. By letting ${\uparrow}x=\bigcap_{V\in \mathsf{ClopUp}(X),x\in V}V$, we can extract a finite clopen upset $A$ such that $x\in A$, and $M\subseteq -A\cup g^{-1}[-U\cup V]$. Hence consider
    \begin{equation*}
S=[-A\cup -g^{-1}[U_{i}]\cup g^{-1}[V_{i}]]\cap \langle A\cap  {\downarrow}(g^{-1}[U_{i}]-g^{-1}[V_{i}])\rangle.
\end{equation*}
Because $X$ has relative pseudocomplements of sets of the form $g^{-1}[U]$, by Proposition \ref{Specific Relative pseudocomplements} we have that ${\downarrow}(g^{-1}[U_{i}]-g^{-1}[V_{i}])$ will be clopen in $X$. Hence the set $S$ is an intersection of subbasic sets. Moreover, $M$ belongs there. Now if $N$ belongs there, then first $N\subseteq -A\cup -g^{-1}[U_{i}]\cup g^{-1}[V_{i}]$, and by assumption, it has some point $x'\in A\cap {\downarrow}g^{-1}[U]-g^{-1}[V]$. Hence $x'\leq y'$ and $y'\in g^{-1}[U]-g^{-1}[V]$. If $x'\leq z'$ and $z'\in N$, then since $x'\in A$, $z'\in A$, so $z'\in g^{-1}[-U\cup V]$. So there can be no point in the same fiber of $g$, i.e., $g(y')\neq g(z')$ whenever $x'\leq z'$ and $z'\in N$. This shows that $M$ has a neighbourhood completely outside of $V_{g}(X)$, showing the latter is closed, as desired.
\end{proof}

We note that the hypothesis of the previous lemma -- that $X$ will have relative pseudocomplements of the form $g^{-1}[U]\rightarrow g^{-1}[V]$ -- will always be satisfied in our contexts. This is the sense in which we will think of $g$ as a way to parametrise those pseudocomplements which we wish to preserve.

\begin{lemma}\label{The root map is good}
The map $r_{g}:V_{g}(X)\to X$ which sends a rooted, closed $g$-open subset to its root, is a continuous, order-preserving and surjective $g$-open map.
\end{lemma}
\begin{proof}
Simply note that if $U$ is a clopen upset, $r_{g}^{-1}[U]=\{M : M\subseteq U\}=[U]$, and $r_{g}^{-1}[X-U]=\{M : M\cap X-U\neq \emptyset\}=\langle X-U\rangle$. The order-preservation is down to the order being reverse inclusion, and the surjectivity follows because, for each $x\in X$, ${\uparrow}x$ will always be $g$-open.

Moreover, note that $r_{g}$ will be $g$-open: if $r_{g}(C)\leq y$, then $m\in C$ is such that $m\leq y$; so because $C$ is $g$-open, there is some $m'$ such that $m\leq m'$ and $g(m')=g(y)$. But then $C'\coloneqq {\uparrow}m'\cap C$ is such that $C\preceq C'$ and $g(r_{g}(C'))=g(y)$, as desired.\end{proof}

\subsection{Free Heyting Algebras from Distributive Lattices}

We will now put the tools developed in the previous section together:

\begin{definition}\label{Ghilardi complex definition}
    Let $g:X\to Y$ be a monotone and continuous function between Priestley spaces $X$ and $Y$. The \textit{$g$-Ghilardi complex} over $X$ $(V_{\bullet}^{g}(X),\leq_{\bullet})$, is an infinite sequence
    \begin{equation*}
        (V_{0}(X),V_{1}(X),...,V_{n}(X),...)
    \end{equation*}
    connected by morphisms $r_{i}:V_{i+1}(X)\to V_{i}(X)$ such that:
    \begin{enumerate}
        \item $V_{0}(X)=X$;
        \item $r_{0}=g$
        \item For $i\geq 0$, $V_{i+1}(X)\coloneqq V_{r_{i}}(V_{i}(X))$;
        \item $r_{i+1}=r_{r_{i}}:V_{i+1}(X)\to V_{i}(X)$ is the root map.
    \end{enumerate}
\end{definition}

Given a $g$-Ghilardi complex over $X$, one can form the projective limit of the sequence, which, since Priestley spaces are closed under such a construction, will again be a Priestley space. We denote this limit by $V_{G}^{g}(X)$ (the $G$ standing for ``Ghilardi"). When $g$ is the terminal map, sending $X$ to $\{\bullet\}$, we drop the superscript, denoting it simply as $V_{G}(X)$. The purpose of such a construction lies in the universal property which is carried out over any one step, which we proceed to outline:

\begin{lemma}\label{Key Lemma in Universal Construction}
Let $g:X\to Y$ be a continuous, order-preserving map. Then the pair $\langle \mathsf{ClopUp}(V_{g}(X)),r_{g}^{-1}\rangle$ has the following universal property: suppose we are given any other pair
\begin{equation*}
\langle D, \mu:\mathsf{ClopUp}(X)\to D\rangle
\end{equation*}
such that
\begin{enumerate}
\item $D$ is a distributive lattice containing relative pseudocomplements of the kind $\mu(C_{1})\rightarrow \mu(C_{2})$ for $C_{1},C_{2}\in \mathsf{ClopUp}(X)$.
\item $\mu(g^{-1}[D_{1}]\rightarrow g^{-1}[D_{2}])=\mu(g^{-1}[D_{1}])\rightarrow \mu(g^{-1}[D_{2}])$ for all $D_{1},D_{2}\in \mathsf{ClopUp}(Y)$.
\end{enumerate}
Then there exists a unique lattice homomorphism $\mu':\mathsf{ClopUp}(V_{g}(X))\to D$ such that the triangle in Figure \ref{fig:commutingtriangleofdistributivelattices} commutes, and such that $\mu'(r^{-1}[C_{1}]\rightarrow r^{-1}[C_{2}])=\mu(C_{1})\rightarrow \mu(C_{2})=\mu'(r^{-1}[C_{1}])\rightarrow \mu'(r^{-1}[C_{2}])$ for all $C_{1},C_{2}\in \mathsf{ClopUp}(X)$.

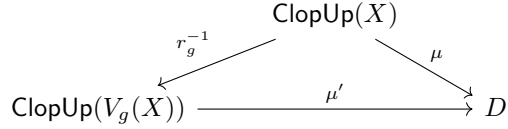
\begin{figure}[h]
\centering
\begin{tikzcd}
& \mathsf{ClopUp}(X) \arrow[rd, "\mu"] \arrow[ld, "r_{g}^{-1}"'] &   \\
\mathsf{ClopUp}(V_{g}(X)) \arrow[rr, "\mu'"] &                                                            & D
\end{tikzcd}
\caption{Commuting Triangle of Distributive Lattices}
\label{fig:commutingtriangleofdistributivelattices}
\end{figure}
\end{lemma}

The proof of this lemma is given by first realising what is the appropriate dual statement.

\begin{lemma}\label{Duality Lemma for Key Property}
The property of Lemma \ref{Key Lemma in Universal Construction} is equivalent to the following: given a Priestley space $Z$ with a $g$-open continuous and order-preserving map $h:Z\to X$, there exists a unique $r_{g}$-open, continuous and order-preserving map $h'$ such that the triangle in Figure \ref{fig:commutingtriangleofpriestleyspaces} commutes.

\begin{figure}[h]
\centering
\begin{tikzcd}
Z \arrow[rd, "h"'] \arrow[rr, "h'"] &   & V_{g}(X) \arrow[ld, "r"] \\                    & X &                         
\end{tikzcd}
\caption{Commuting Triangle of Priestley spaces}
\label{fig:commutingtriangleofpriestleyspaces}
\end{figure}
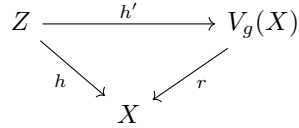
\end{lemma}
\begin{proof}
This follows immediately from Lemma \ref{G-open equivalent to preservation of pseudocomplements} by dualising the statements.
\end{proof}

\begin{proof}(Of Lemma \ref{Key Lemma in Universal Construction})
We show the property from Lemma \ref{Duality Lemma for Key Property}. Given $Z$, $X$ and $h$ as given, first we note that the definition of $h'$ is totally forced by the commutativity of the triangle: if $k:Z\to V_{g}(X)$ is any map in these conditions, we will show that
\begin{equation*}
    k(a)=h[{\uparrow}a].
\end{equation*}
Indeed, since the diagram commutes, the root of $k(a)$ will be $h(a)$; and if $a\leq b$, then $k(a)\leq k(b)$, so $h[{\uparrow}a]\subseteq k(a)$. Now if $x\in k(a)$, consider $k(a)\cap {\uparrow}x$. Then $k(a)\leq k(a)\cap {\uparrow}x$, so because $k$ is $r_{g}$-open, there is some $a'\geq a$ such that $r_{g}(k(a'))=r_{g}(k(a)\cap {\uparrow}x)=x$. Hence $h(a')=x$, and so $x=h(a')\in h[{\uparrow}a]$, as desired.

It thus suffices to show that given $a\in Z$, the assignment
\begin{equation*}
h'(a)=\{h(b) : a\leq b\},
\end{equation*}
indeed defines a function as desired. Note that $h'(a)=h[{\uparrow}a]$, and since $h$ is continuous between Stone spaces, and hence a closed map, $h'$ maps $a$ to a closed subset; it is of course also rooted and $h(a)$. And it is $g$-open as a subset, because $h$ is $g$-open as a map. Clearly $r_{g}(h'(a))=h(a)$. It is continuous, since $(h')^{-1}[[-U\cup V]]=\{a : h'(a)\in [-U\cup V]\}$, and this is the same as saying that ${\uparrow}a\subseteq h^{-1}[-U]\cup h^{-1}[V]$. Now since $Z$ has relative pseudocomplements of these sets, then $\Box(-U\cup V)$ exists, and indeed
\begin{equation*}
\{a : {\uparrow}a\subseteq h^{-1}[-U]\cup h^{-1}[V]\}=\{a : a\in \Box(h^{-1}[-U]\cup h^{-1}[V])\}=\Box(h^{-1}[-U]\cup h^{-1}[V]).
\end{equation*}
This shows continuity. It is also of course order-preserving, since if $a\leq b$, then ${\uparrow}b\subseteq {\uparrow}a$, so $h[{\uparrow}b]\subseteq h[{\uparrow}a]$, i.e., $h'(a)\leq h'(b)$. Finally it is also $r_{g}$-open: if $a\in Z$, and $h'(a)\leq M$ where $M\in V_{g}(X)$, then $M\subseteq h'(a)$, so if $z$ is the root of $M$, then $z=h(c)$, where $a\leq c$. Then $r(h'(c))=r(M)$, as desired.\end{proof}

Moreover, we can see that the construction $V_{g}$ also works functorially:

\begin{lemma}\label{Lifting on functions}
    Let $p:X\to Y$ be a monotone and continuous function between Priestley spaces, and $g_{X}:X\to Z$ and $g_{Y}:Y\to Z$ be two maps such that the obvious triangle commutes, and ${\downarrow}h^{-1}[U-V]$ is clopen whenever $U,V$ are clopen upsets and $h\in \{g_{X},g_{Y}\}$. Then $p[-]:V_{g}(X)\to V_{g}(Y)$ is a unique monotone and continuous function such that for each $C\in V_{g_{x}}(X)$, $pr(C)=r_{g_{Y}}p[C]$.
\end{lemma}
\begin{proof}
    It is clear that $p[-]$ will be a monotone and continuous function. To see its uniqueness, note that $p[-]$ is $r_{g_{Y}}$-open, and hence by Lemma \ref{Key Lemma in Universal Construction}, it is the unique extension of $pr$ (which is $g_{Y}$-open given that $p$ is $g_{Y}$-open). 
\end{proof}

Recall that since $V_{G}(X)$ is a projective limit, there is a map $\pi_{0}:V_{G}(X)\to X$ which projects an inverse limit to $X$. Note that $\pi_{0}$ is actually surjective, since for each $x\in X$, we can consider the sequence:
\begin{equation*}
    (x,{\uparrow}x,{\uparrow}({\uparrow}x),...)
\end{equation*}
which will be an element in $V_{G}(X)$, mapping to $x$.

\begin{lemma}\label{Lifting Lemma on the Esakia spaces}
    Let $g:X\to Y$ be a monotone and continuous function such that $g$-indexed relative pseudocomplements exist. Suppose that $Z$ is an Esakia space and $h:Z\to X$ is a monotone and continuous function which is $g$-open. Then there exists a unique p-morphism $\overline{h}:Z\to V_{G}(X)$ such that $\pi_{0}\overline{h}=h$. 
\end{lemma}
\begin{proof}
    Using repeatedly Lemma \ref{Key Lemma in Universal Construction}, given that $Z$ is an Esakia space, we get a sequence of maps, starting with $h_{0}=h$, and letting $h_{n+1}=h_{n}[-]$. Then for each $x\in Z$, we have that $\overline{h}(x)$ will be the element
    \begin{equation*}
        (h_{0}(x),h_{1}(x),...,h_{n}(x),...).
    \end{equation*}
    It is clear that this map is an order-preserving map and continuous, by the universal property of inverse limits. We check that it is a p-morphism. Assume that for $x\in Z$, $\overline{h}(x)\leq y$. Then consider the following set:
    \begin{equation*}
        {\uparrow}x\cap \bigcap \{h_{n}^{-1}[y(n)] : n\in \omega\}.
    \end{equation*}
    Since $Z$ is compact, and all the sets involved are closed, this intersection will be empty if and only if a finite intersection is empty. Let $n$ be such that
    \begin{equation*}
        {\uparrow}x\cap h_{0}^{-1}[y(0)]\cap...\cap h_{n}^{-1}[y(n)]=\emptyset.
    \end{equation*}
    Since $\overline{h}(x)\leq y$, then note that $h_{n+1}(x)\leq y(n+1)$; so there is some $k\geq x$ such that $h_{n}(k)=y(n)$. Then certainly $k\in {\uparrow}x\cap h_{0}^{-1}[y(0)]\cap...\cap h_{n}^{-1}[y(n)]$, which contradicts the above set being empty. By reductio we get $x\leq k$ such that $\overline{h}(k)=y$, as desired.
\end{proof}

Finally, we can lift monotone and continuous functions satisfying appropriate $g$-openness conditions functorially:

\begin{lemma}\label{Lemma on lifting maps}
    Let $p:X\to Y$ be a monotone and continuous function between Priestley spaces which is $g_{Y}$-open, where $g_{X}:X\to Z$ and $g_{Y}:Y\to Z$ are two maps such that the obvious triangle commutes and the relative pseudocomplements indexed by $g_{X}$ and $g_{Y}$ exist. Then there is a unique p-morphism $\overline{p}:V_{G}(X)\to V_{G}(Y)$ such that $\pi_{0}^{Y}\overline{p}=p\pi_{0}^{X}$.
\end{lemma}
\begin{proof}
    We use Lemma \ref{Lifting on functions} repeatedly on the $g_{X}$-Ghilardi complex over $X$ and the $g_{Y}$-Ghilardi complex over $Y$: we let $p_{0}=p$, and $p_{n+1}=p[-]$. This lifts to a function on the inverse limits mapping $x\in V_{G}(X)$ to
    \begin{equation*}
        (p_{0}(x(0)),p_{1}(x(1)),...,p_{n}(x(n)),...).
    \end{equation*}
    It is not difficult to see that having $\overline{p}$ so defined, $\pi_{0}\overline{p}$ is a $g_{Y}$-open map, so $\overline{p}$ will additionally coincide with the unique p-morphism lifting of $p\pi_{0}$ (by Lemma \ref{Lifting Lemma on the Esakia spaces}), in addition to making all diagrams commute. 
\end{proof}

We can now define $V_{G}:\mathbf{Pries}\to \mathbf{Esa}$ as a functor, assigning $V_{G}(X)$ to each Priestley space $X$, and $\overline{p}:V_{G}(X)\to V_{G}(Y)$ to each monotone and continuous function $p:X\to Y$. We then obtain:

\begin{theorem}\label{Adjunction Theorem}
    The functor $V_{G}$ is right adjoint to the inclusion $I:\mathbf{Esa}\to \mathbf{Pries}$. 
\end{theorem}
\begin{proof}
    This follows immediately from Lemma \ref{Lifting Lemma on the Esakia spaces} and Lemma \ref{Lemma on lifting maps}. Note that the unit $\eta_{X}:X\to V_{G}(X)$ consists of the unique lifting of the identity $id_{X}$, and is an embedding in the category of Esakia spaces; the co-unit $\varepsilon_{X}:V_{G}(X)\to X$ is the projection onto the first coordinate, and is surjective as noted before. Hence the right adjoint is full whilst the left adjoint is faithful.
\end{proof}

Note that the inclusion of $\mathbf{Esa}$ into $\mathbf{Pries}$ is certainly not full, since there are monotone maps between Esakia spaces which are not p-morphisms. And the co-unit can certainly fail to be an isomorphism, since given $\mathbf{2}$, the poset $0<1$, $V_{G}(\mathbf{2})$ coincides with the Rieger-Nishimura ladder, which is infinite.

We can make sense of the above algebraically as well. Let $\alg{D}$ be a distributive lattice. Then we can construct $\alg{D}_{0}=\alg{D}$ and $\alg{D}_{n+1}$ as specified in subsection \ref{Conceptual idea}: we form
\begin{equation*}
    \mathbf{F}_{DLat}(\{(a,b) : a,b\in D_{n}\})/\Theta
\end{equation*}
i.e., the free distributive lattice on the pairs $a,b\in D_{n}$, quotiented by the congruence forcing these pairs to act like relative pseudocomplements, and forcing the relative pseudocomplements added at the stage $n+2$ to agree with those added at stage $n+1$. We can then consider the homomorphism $i_{n}:D_{n}\to D_{n+1}$ which sends $a$ to the equivalence class $[(1,a)]_{\Theta}$. This gives us a directed system $(D_{n},i_{n})$ which, as is easy to see, has the same universal properties as the Ghilardi complex over $X$, where $X$ is the dual Priestley space of $D$, and in fact this means that $i_{n}$ are all inclusions. By Lemma \ref{Directed Unions are dual to inverse limits}, we have that forming the directed limit is dual to taking the inverse limit, i.e., $\lim_{n\in \omega}D_{n}$ is dual to $V_{G}(X)$. Let $F_{G}(D)=\lim_{n\in \omega}D_{n}$ as stated. Moreover, $F_{G}$ forms a functor as well, obtained by dualizing all the above facts. Hence we have:

\begin{theorem}\label{Key Theorem on free Heyting algebras}
The functor $F_{G}:\mathbf{DLat}\to \mathbf{HA}$ is left adjoint to the inclusion of Heyting algebras in Distributive lattices.
\end{theorem}



    

We finish with an explicit calculation of $V_{G}$ that illustrates the general case:

\begin{example}
    Let $X$ be the Priestley space in Figure \ref{fig:spaceX}.

    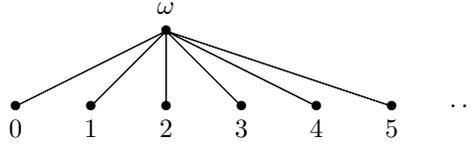
\begin{figure}
        \centering
\begin{tikzpicture}
    \node at (0,0) {$\bullet$};
    \node at (0,-0.3) {$2$};
    \node at (0,1) {$\bullet$};
    \node at (-1,0) {$\bullet$};
    \node at (-1,-0.3) {$1$};
    \node at (-2,0) {$\bullet$};
    \node at (-2,-0.3) {$0$};
    \node at (1,0) {$\bullet$};
    \node at (0,1.3) {$\omega$};
    \node at (1,-0.3) {$3$};
    \node at (2,0) {$\bullet$};
    \node at (2,-0.3) {$4$};
    \node at (3,0) {$\bullet$};
    \node at (3,-0.3) {$5$};
    \node at (4,0) {$\dots$};

    \draw (0,1) -- (-2,0) -- (0,1) -- (-1,0) -- (0,1) -- (0,0) -- (0,1) -- (1,0) -- (0,1) -- (2,0) -- (0,1) -- (3,0) -- (0,1);
\end{tikzpicture}        \caption{The space $X$}
        \label{fig:spaceX}
    \end{figure}
    We endow this space with a Priestley topology by declaring a basis consisting of the finite sets and the cofinite sets containing $\omega$. Now consider $V_{1}(X)$. Closed and rooted subsets will be precisely the singletons of the natural numbers, together, with $\{\omega\}$, and the pairs of the form $\{n,\omega\}$. Hence $V_{1}(X)$ will have the poset structure as detailed in Figure \ref{fig:spaceV1X}

        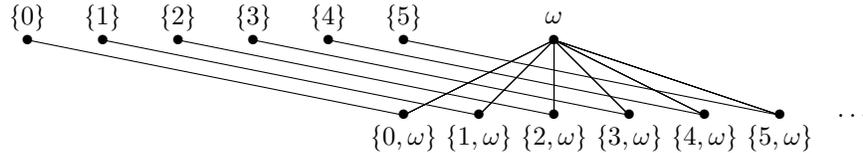
\begin{figure}[h]
        \centering
\begin{tikzpicture}
    \node at (0,0) {$\bullet$};
    \node at (0,-0.3) {$\{2,\omega\}$};
    \node at (0,1) {$\bullet$};
    \node at (-1,0) {$\bullet$};
    \node at (-1,-0.3) {$\{1,\omega\}$};
    \node at (-2,0) {$\bullet$};
    \node at (-2,-0.3) {$\{0,\omega\}$};
    \node at (1,0) {$\bullet$};
    \node at (0,1.3) {$\omega$};
    \node at (1,-0.3) {$\{3,\omega\}$};
    \node at (2,0) {$\bullet$};
    \node at (2,-0.3) {$\{4,\omega\}$};
    \node at (3,0) {$\bullet$};
    \node at (3,-0.3) {$\{5,\omega\}$};
    \node at (4,0) {$\dots$};

    \draw (0,1) -- (-2,0) -- (0,1) -- (-1,0) -- (0,1) -- (0,0) -- (0,1) -- (1,0) -- (0,1) -- (2,0) -- (0,1) -- (3,0) -- (0,1);

    \node at (-7,1) {$\bullet$};
    \draw (-2,0) -- (-7,1);
    \node at (-7,1.3) {$\{0\}$};
    \node at (-6,1) {$\bullet$};
    \draw (-1,0) -- (-6,1);
    \node at (-6,1.3) {$\{1\}$};
    \node at (-5,1) {$\bullet$};
    \draw (0,0) -- (-5,1);
    \node at (-5,1.3) {$\{2\}$};
    \node at (-4,1) {$\bullet$};
    \draw (1,0) -- (-4,1);
    \node at (-4,1.3) {$\{3\}$};
    \node at (-3,1) {$\bullet$};
    \draw (2,0) -- (-3,1);
    \node at (-3,1.3) {$\{4\}$};
    \node at (-2,1) {$\bullet$};
    \draw (3,0) -- (-2,1);
    \node at (-2,1.3) {$\{5\}$};
    
\end{tikzpicture}        \caption{The space $V_{1}(X)$}
        \label{fig:spaceV1X}
    \end{figure}

    The topology will again be given by the finite subsets not containing $\omega$, and the cofinite subsets containing it. Proceeding in this way, we obtain spaces $V_{2}(X)$, $V_{3}(X)$, etc, such that each of them is topologically isomorphic to the Alexandroff compactification of $\omega$, and order theoretically, on each root, it looks like the $n$-th step construction of the Rieger-Nishimura. Then $V_{G}(X)$ will be isomorphic to infinitely many copies of the Rieger-Nishimura, where the point labelled with $1$ is identified between all of them, and compactifies the disjoint union of these topologies.
\end{example}

We will have opportunity in the subsequent sections to see further examples of computations involving this left adjoint, including some refinements that will be necessary moving forward. For now we will turn to some basic applications of these constructions.

\section{Some Applications to the theory of Heyting Algebras }\label{Applications}

In this section we provide some basic applications of the construction of Section \ref{Free Constructions of Heyting Algebras over Distributive Lattices} to the theory of Heyting algebras. Whilst in subsequent sections we will turn to some novel applications, in this section we will derived known results (or results easily derived from results in the literature); however, we believe that the methodology can be suggestive to tackle similar questions in settings which may not have such a developed theory as Heyting algebras.  

\subsection{Free Heyting Algebras}

As a first application of Theorem \ref{Key Theorem on free Heyting algebras}, we show how to obtain a description of the free Heyting algebra on any number of generators, and study some of its basic properties. Let $X$ be an arbitrary set, and recall Definition \ref{Description of free algebras}, i.e., $\mathbf{D}(X)$ being the free distributive lattice on $X$ many generators.

\begin{theorem}
    The algebra $F_{G}(\mathbf{D}(X))$ is the free Heyting algebra on $X$ many generators.
\end{theorem}
\begin{proof}
    Since $F_{G}(\mathbf{D}(-))$ is the composition of the left adjoint from the inclusion of $\mathbf{DLat}$ in $\mathbf{Set}$, and $\mathbf{HA}$ in $\mathbf{DLat}$, it will be left adjoint to the composed inclusion. But by uniqueness of adjoints, it will be the free functor on the forgetful functor from $\mathbf{HA}$ to $\mathbf{DLat}$.
\end{proof}

This construction restricts, in the case of $X$ being finite, to the Ghilardi/Urquhart construction \cite{ghilardifreeheyting,Urquhart1973}. We note that unlike that situation, it is not immediately clear that $F_{G}(\mathbf{D}(X))$ carries the additional structure of a \textit{bi-Heyting algebra}; that is because, to show it, it would be necessary to show that $V_{1}(2^{X})$ is likewise a bi-Esakia space. A natural associated question is whether any of the free Heyting algebras is complete. We can use some well-known auxiliary constructions to show that this is never the case for infinitely many $X$:

\begin{definition}
    Let $X$ be an Esakia space. We say that $X$ is \textit{extremally order-disconnected} if whenever $U$ is an open upset, then $cl(U)$ is clopen.
\end{definition}

It was shown in \cite{Bezhanishvili2008} that a Heyting algebra $\alg{H}$ is complete if and only if its dual space $X$ is extremally order-disconnected. Given an ordered-topological space $X$, we write $\max(X)$ for the maximum of $X$. If $X$ is a Priestley space, it is well-known \cite{Esakiach2019HeyAlg} that $\max(X)$ is a Stone space.

\begin{lemma}\label{Cutting down an extremally order-disconnected yields an extremally disconnected}
    If $X$ is extremally order-disconnected, then $\max(X)$ is extremally disconnected.
\end{lemma}
\begin{proof}
    Suppose that $U\subseteq \max(U)$ is an open set. Then by definition, there is some open subset $U'$ such that $U=U'\cap \max(X)$. Look at $\Box U'$, which is open, since $X$ is Priestley. Then note that $\Box U'\cap \max(X)=U$: indeed, one inclusion is clear since $\Box U'\subseteq U'$, and for the other, if $x\in U$, then $x\in U'$, and certainly whenever $x\leq y$, $y\in U'$, since $x$ is maximal. Hence we can take $U'$ to be an open upset.

    We wish to show that $cl(U)$ is clopen. Note that $cl(U)=cl(U'\cap \max(X))=cl(U')\cap \max(X)$, since we are sitting in the maximum. So since $U'$ is an open upset, $cl(U')$ is clopen, and $\max(X)-cl(U)=\max(X)-cl(U')$. This shows the result.    
\end{proof}

We will need the following lemma:

\begin{lemma}\label{Homeomorphisms of Esakia spaces and their maxima}
    Let $X$ be an Esakia space. Then $\max(V_{G}(X))$ is homeomorphic to $X$.
\end{lemma}
\begin{proof}
    Consider the map $i:X\to V_{1}(X)$ defined by sending $x$ to $\{x\}$. Note that this is continuous: for each clopen upset $U$ in $X$, $\{x\}\in [U]$ if and only if $\{x\}\in \langle U\rangle$ if and only if $x\in U$. It is trivially bijective, hence, it is a homeomorphism, since $\max(V_{1}(X))$ is a Stone space. Now we prove by induction that that $x\in \max(V_{n+1}(X))$ if and only if $\{x\}\in \max(V_{n+2}(X))$, and $r_{n}(x)\in \max(V_{n}(X))$: if $x\in \max(V_{n+1}(X))$, then ${\uparrow}x=\{x\}$ is a rooted and $r_{n}$-open subset; conversely, if $\{x\}\in \max(V_{n+1}(X))$, and $x\leq y$, then by $r_{n}$-openness, there ought to be some $k$, such that $x\leq k$, $k\in \{x\}$ and $r_{n}(k)=r_{n}(y)$. Hence $r_{n}(x)=r_{n}(y)$ which is a maximal element, i.e., $x=y$. Now if we consider $\max(V_{G}(X))$, note that $x\in \max(V_{G}(X))$ if and only if $x(i)\in \max(V_{i}(X))$; hence this will again consist of the sequences:
    \begin{equation*}
        (x,\{x\},\{\{x\}\},...)
    \end{equation*}
    for $x\in X$. Thus $\max(V_{G}(X))\cong \max(V_{1}(X))\cong X$, as desired.
\end{proof}

\begin{corollary}
    If $|X|>1$, $F_{G}(\mathbf{D}(X))$ is not complete.
\end{corollary}
\begin{proof}
For $X$ finite, this was shown by Bellissima \cite[Theorem 4.2]{Bellissima1986}. We show this for infinite $X$. 

By Lemma \ref{Homeomorphisms of Esakia spaces and their maxima}, $\max(V_{G}(\mathbf{2}^{X}))=\mathbf{2}^{X}$, i.e., the maximum is homeomorphic to the generalized Cantor space on $X$ many elements. Were $V_{G}(\mathbf{2}^{X})$ extremally order-disconnected, by Lemma \ref{Cutting down an extremally order-disconnected yields an extremally disconnected} we would have that $\mathbf{2}^{X}$ would be extremally disconnected, which is not true: dually, the free Boolean algebras on $X$ generators would be complete Boolean algebras, which is known to not hold \cite[pp.133]{Koppelberg1989-lg}.
\end{proof}

\subsection{Coproducts and Colimits of Heyting Algebras}

As a further application, we can provide an explicit description of the coproduct of two Heyting algebras.

For that purpose, let $X,Y$ be two Esakia spaces. Note that $X\times Y$, given the product topology and order, is again an Esakia space \cite{Esakiach2019HeyAlg}, and in fact the maps $\pi_{X}:X\times Y\to X$ and $\pi_{Y}:X\times Y\to Y$ are p-morphisms. Nevertheless it is not difficult to see that $X\times Y$ may fail to be the categorical product:

\begin{example}
    Consider the figure \ref{fig:Product of Esakia spaces which is not their categorical product}.

\begin{figure}[h]
\centering
\begin{tikzcd}
            &             &                    & {(x,a)}                       &                    \\
x           & a           & {(y,a)} \arrow[ru] &                               & {(x,b)} \arrow[lu] \\
y \arrow[u] & b \arrow[u] &                    & {(y,b)} \arrow[lu] \arrow[ru] &                   
\end{tikzcd}
\caption{Product of Esakia spaces which is not the categorical product}
\label{fig:Product of Esakia spaces which is not their categorical product}
\end{figure}
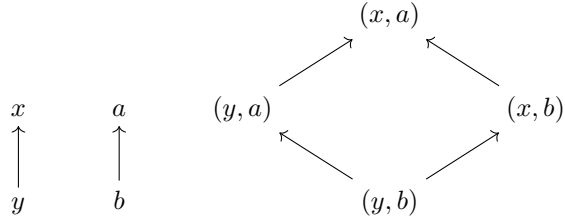

Let $P_{1}=\{x,y\}$ and $P_{2}=\{a,b\}$. To see that $P_{1}\times P_{2}$ is not the categorical product, consider $P_{1}\sqcup P_{2}$, the disjoint union of the two posets, and the following two maps: $f:P_{1}\sqcup P_{2}\to P_{2}$ sends $P_{1}$ to $a$, and maps $P_{2}$ as the identity; and $g:P_{1}\sqcup P_{2}\to P_{1}$ maps $P_{1}$ as the identity, and sends $b$ to $y$ and $a$ to $x$. Hence $f(b)=b$ and $g(b)=y$, so in the universal map $(g\times f)(b)=(y,b)$. Similarly, $(g\times f)(y)=(y,a)$, $(g\times f)(a)=(x,a)$ and $(g\times f)(b)=(y,b)$. Then $k(b)\leq (y,a)$ but whenever $b\leq z$, then $k(z)\neq (y,a)$. So $k$ is not a p-morphism.
\end{example}

The key problem is that $X\times Y$ is not sufficiently general to be able to do the job of the product. If $X=V_{G}(X')$ and $Y=V_{G}(Y')$ we would know how to solve the problem -- using the preservation of products by right adjoints, we would have that $X\otimes Y \cong V_{G}(X\times Y)$. Whilst this idea is in general too coarse, with a bit more care, we can obtain the correct description of the categorical product of Esakia spaces, which dually gives us the coproduct of Heyting algebras.

For that purpose we will need Remark \ref{Remark on multiple maps}: given $X,Y$ we will consider the space
\begin{equation*}
    V_{1}^{\pi_{X},\pi_{Y}}(X\times Y)
\end{equation*}
consisting of the rooted, closed  $\pi_{X}$-open and $\pi_{Y}$-open subsets of $X\times Y$. It follows from Lemmas \ref{The root map is good} and \ref{The root map is good} that $V_{1}^{\pi_{X},\pi_{Y}}(X\times Y)$ is a Priestley space, and that the root map is continuous, monotone and $\pi_{X},\pi_{Y}$-open.

We will then define the $(\pi_{X},\pi_{Y})$-Ghilardi complex over $X\times Y$ in an analogous way to Definition \ref{Ghilardi complex definition}, with the only difference being that we have two maps at the base. Let $V_{G}^{\pi_{X},\pi_{Y}}(X\times Y)$ be the inverse limit of that inverse system. Then we have the following:

\begin{lemma}
    The maps $\pi_{X}\pi_{0}:V_{G}^{\pi_{X},\pi_{Y}}(X\times Y)\to X$ and $\pi_{Y}\pi_{0}:V_{G}^{\pi_{X},\pi_{Y}}(X\times Y)\to Y$ are monotone and continuous p-morphisms.
\end{lemma}
\begin{proof}
    It is clear that they are continuous and monotone; we focus on showing the p-morphism condition for $X$. Assume that $C\in V_{G}^{\pi_{X},\pi_{Y}}(X\times Y)$, and $\pi_{X}\pi_{0}(C)\leq z$. Note that $C(0)$ is of the form $(x,y)$ for $x\in X$ and $y\in Y$. Since $C(1)$ is $\pi_{X}$-open, and $(x,y)\in C(1)$, and $(x,y)\leq (z,y)$, there must be some $(z,y')\in C(1)$ such that $(x,y)\leq (z,y')$.

    Similarly, since $(z,y')\in C(1)$, we can find $K=C(1)\cap {\uparrow}(z,y')$, and have that $C(1)\leq K$; since $C(2)$ is $r_{1}$-open, there is some $K'\in C(2)$ such that $C(1)\leq K'$ and $r_{1}(K')=r_{1}(K)=(z,y')$. In general, letting $K_{0}=(z,y')$ and $K_{1}=K'$, we can use the fact that $K_{n}\in C(n+1)$ to obtain some $K_{n+1}\in C(n+2)$ in the same conditions. In this way, we obtain a chain
    \begin{equation*}
        Z=(K_{0},K_{1},K_{2},...,K_{n},...),
    \end{equation*}
    with the property that $C\leq Z$, and such that $\pi_{X}\pi_{0}(Z)=z$. This shows that $\pi_{X}\pi_{0}$ is a p-morphism as desired.
\end{proof}

Let $p_{X}=\pi_{X}\pi_{0}$ and $p_{Y}=\pi_{Y}\pi_{0}$. We can thus consider the triple $(V_{G}^{\pi_{X},\pi_{Y}}(X\times Y),p_{X},p_{Y})$ and have the following:

\begin{theorem}\label{Characterization of the categorical product}
    Given two Esakia spaces $X$ and $Y$, the triple $(V_{G}^{\pi_{X},\pi_{Y}}(X\times Y),p_{X},p_{Y})$ is the categorical product of $X$ and $Y$.
\end{theorem}
\begin{proof}
    Given any Esakia space $Z$, if $f:Z\to X$ and $g:Z\to Y$ are two p-morphisms, then $(f\times g):Z\to X\times Y$ is a $\pi_{X},\pi_{Y}$-open map, and $\pi_{X}(f\times g)=f$ and $\pi_{Y}(f\times g)=g$. Using the universal property of $V_{n}^{\pi_{X},\pi_{Y}}(X\times Y)$, we obtain a unique map to each such space, such that $(f\times g)_{0}=(f\times g)$ and $(f\times g)_{n+1}=(f\times g)_{n}[-]$; this in turn forces, by Lemma \ref{Duality Lemma for Key Property} and the properties of inverse limits,  the existence of a unique p-morphism to $(f\times g)_{\infty}:Z\to V_{G}^{\pi_{X},\pi_{Y}}(X\times Y)$, commuting with this diagram, i.e., such that $\pi_{0}(f\times g)_{\infty}=(f\times g)$. Then $p_{X}(f\times g)_{\infty}=f$ and $p_{Y}(f\times g)_{\infty}=g$, as desired.
\end{proof}

The reader will observe that, in line with our comments, the former is not limited to binary products, which means that we can provide an explicit description of arbitrary products of Esakia spaces as well. More interestingly, having such an explicit description of the coproduct allows us to prove some purely categorical properties of the category of Heyting algebras.  

\begin{definition}
    Let $\mathbf{C}$ be a category with finite products ($\times$) and coproducts $(\oplus)$. We say that $\mathbf{C}$ is co-distributive if coproducts distribute over products, i.e., for $\alg{A},\alg{B},\alg{C}\in \mathbf{C}$ we have
    \begin{equation*}
        \alg{A}\oplus (\alg{B}\otimes \alg{C})\cong (\alg{A}\oplus \alg{B})\times (\alg{A}\oplus \alg{C}).
    \end{equation*}
\end{definition}

\begin{theorem}
The category of Heyting algebras is co-distributive.
\end{theorem}
\begin{proof}
Dually denote by $X\otimes Y$ the categorical product of two Esakia spaces. Let $X,Y,Z$ be three Esakia spaces, we will show that
\begin{equation*}
    X\otimes (Y\sqcup Z)\cong (X\otimes Y)\sqcup (X\otimes Z).
\end{equation*}

First note that taking the product in Priestley spaces, we have $X\times (Y\sqcup Z)\cong (X\times Y)\sqcup (X\times Z)$, a fact which follows easily from the corresponding fact for sets. Now we will show that for each $n$:
\begin{equation*}
    V_{n}^{\pi_{X},\pi_{Y\sqcup Z}}(X\times (Y\sqcup Z))\cong V_{n}^{\pi_{X},\pi_{Y}}(X\times Y)\sqcup V_{n}^{\pi_{X},\pi_{Z}}(X\times Z).
\end{equation*}
For $n=0$ this is the aforementioned fact. Now assume that it holds for $n$. Given $C\in V_{n+1}^{\pi_{X},\pi_{Y\sqcup Z}}(X\times (Y\sqcup Z))$, note that $C$ is a subset of $V_{n}^{\pi_{X},\pi_{Y\sqcup Z}}(X\times (Y\sqcup Z))$; the latter is by induction hypothesis isomorphic to $V_{n}^{\pi_{X},\pi_{Y}}(X\times Y)\sqcup V_{n}^{\pi_{X},\pi_{Z}}(X\times Z)$, so via the isomorphism, $C$ is a subset of the disjoint union; but since $C$ is rooted, it is a subset of one of the two disjoints. This assignment is easily seen to be bijective, continuous, and monotone in both directions, giving us an isomorphism.

Moreover, note that inverse limits commute with disjoint unions: $$\varprojlim_{n\in \omega}V_{n}(X\times Y)+V_{n}(X\times Z)\cong (\varprojlim_{n\in \omega}V_{n}(X\times Y))+(\varprojlim_{n\in \omega}V_{n}(X\times X_{C}).$$

To see this, note that if $x$ belongs to the first set, and $x(0)\in V_{0}(X\times Y)$, then since $x(1)$ must have its root in this set, $x(1)\in V_{1}(X\times Y)$; this remains true for each $n$, so this means that $x\in \varprojlim_{n\in \omega}V_{n}(X\times Y)$. Conversely, if $x$ belongs to one of the two sets, then certainly it belongs to the projective limit of the disjoint unions.

Putting these facts together we have that:
\begin{align*}
    X\otimes (Y\sqcup Z) &\cong V_{G}^{\pi_{X},\pi_{Y\sqcup Z}}(X\times (Y\sqcup Z))\\
    &\cong \varprojlim_{n\in \omega}V_{n}^{\pi_{X},\pi_{Y\sqcup Z}}(X\times (Y\sqcup Z))\\
    &\cong \varprojlim_{n\in \omega}V_{n}^{\pi_{X},\pi_{Y}}(X\times Y)\sqcup V_{n}^{\pi_{X},\pi_{Z}}(X\times Z)\\
    &\cong (\varprojlim_{n\in \omega}V_{n}(X\times Y))+(\varprojlim_{n\in \omega}V_{n}(X\times X_{C})\\
    &\cong X\otimes Y\sqcup X\otimes Z.
\end{align*}
\end{proof}

Using similar techniques one can show that coproducts of Heyting algebras are disjoint, and that the category of Heyting algebras form a co-extensive category. This could be deduced from known results in \cite{Carboni1993}, though the explicit proof given here might be of independent interest. Finally we turn to a question of amalgamation: indeed, making use of the explicit description of the product given before, we can likewise provide, dually, a description of the pullback of two Esakia spaces:

\begin{theorem}
    Let $X,Y,Z$ be Esakia spaces, and let $f:X\to Z$ and $g:Y\to Z$ be two p-morphisms. Let $X\times_{Z}Y$ be the pullback in the category of Priestley spaces. Then the pullback of this diagram consists of
    \begin{equation*}
        X\otimes_{Z}Y\coloneqq V_{G}^{\pi_{X},\pi_{Y}}(X\times_{Z}Y)
    \end{equation*}
\end{theorem}

The proof of this follows exactly the same ideas as that of Theorem \ref{Characterization of the categorical product}; indeed note that the projections from $V_{G}^{\pi_{X},\pi_{Y}}(X\times_{Z}Y)$ to $X$ and $Y$ are always surjective, since for each $x\in X$ (resp. $y\in Y$) ${\uparrow}x$ is a $(\pi_{X},\pi_{Y})$-open closed and rooted subset. Using this pullback we can obtain an explicit proof of the amalgamation property for Heyting algebras, which is ``constructive" in the sense that it is brought about by categorical, rather than model-theoretic or logical, considerations.

\begin{theorem}
    The variety of Heyting algebras has the amalgamation property. 
\end{theorem}
\begin{proof}
    Dually, this amounts to showing that whenever we have a cospan diagram $f:X\to Z$ and $g:Y\to Z$, where both $f$ and $g$ are surjective, then there is some $W$ which is a cone for this diagram, and where the maps to $X$ and $Y$ are again surjective. Note that forming the pullback in the category of Priestley spaces, the projection maps
    \begin{equation*}
        X\xleftarrow{\pi_{X}} X\times_{Z} Y \xrightarrow{\pi_{Y}} Y
    \end{equation*}
    are surjective, since for each $x\in X$, $f(x)=g(y)$ for some $y\in Y$, by surjectivity of $g$, and vice-versa. Certainly then the unique lifting from $V_{G}(X\times_{Z}Y,p_{Y},p_{X})$ to the projections will likewise be surjective, and both will be p-morphisms, showing that amalgamation holds.
\end{proof}

\subsection{Esakia Duality from Priestley Duality}

As a final application we show how Esakia duality for Heyting algebras can be presented coalgebraically on the basis of Priestley by using the previous functor. This recovers the results from \cite{daveycoalgebraheyting}. For that purpose, given a Priestley space $X$, recall the map
\begin{equation*}
    \varepsilon:V_{G}(X)\to X,
\end{equation*}
the counit of the adjunction, which is simply $\pi_{0}$. Moreover, if $X$ is an Esakia space, recall the map
\begin{equation*}
    \eta:X\to V_{G}(X)
\end{equation*}
the unit of the adjunction, which is the unique p-morphism lifting the identity. Explicitly, $\eta$ maps $x\in X$ to the tuple
\begin{equation*}
    (x,{\uparrow}x,{\uparrow}({\uparrow}x),...).
\end{equation*}
We write $\mu$ for the map $V_{G}(\eta)$. We will now see the functor $V_{G}$ as the composition $V_{G}\circ I$, where $I$ is the inclusion, i.e., we will see it as a comonad on the category of Priestley spaces.  In this way, recall that a \textit{coalgebra for this comonad} is a Priestley space $X$, together with a monotone and continuous function $f:X\to V_{G}(X)$ subject to the following diagrams commuting:

\begin{figure}[h]
    \centering
\begin{tikzcd}
X \arrow[r, "f"] \arrow[rd, "1_{X}"'] & V_{G}(X) \arrow[d, "\varepsilon_{X}"'] & X \arrow[d, "f"'] \arrow[r, "f"] & V_{G}(X) \arrow[d, "\mu"] \\
& X                                      & V_{G}(X) \arrow[r, "V_{G}(f)"']  & V_{G}(V_{G}(X))          
\end{tikzcd}    \caption{Coalgebras for the $V_{G}$ comonad}
    \label{fig:coalgebrasforVGcomonad}
\end{figure}

Note that the former diagram implies that $\varepsilon_{X}(f(x))=x$. We will now show the following:

\begin{theorem}\label{comonadic adjunction}
    Let $X$ be a Priestley space. Then the following are equivalent:
    \begin{enumerate}
        \item $X$ is an Esakia space;
        \item $X$ admits a (unique) $V_{G}$-coalgebra structure.
        \item The map ${\uparrow}:X\to V(X)$ is a monotone and continuous function;
        \item The map ${\uparrow}:V(X)\to V(X)$ is a monotone and continuous function.
    \end{enumerate}
\end{theorem}
\begin{proof}
(1) implies (2): if $X$ is an Esakia space, then it admits a free $V_{G}$-coalgebra structure, namely the unit map $\eta:X\to V_{G}(X)$. Moreover, this structure is unique: if $f:X\to V_{G}(X)$ is any such coalgebra, then for each $x\in X$, $f(x)=(a_{0},a_{1},...)$; by the commutation conditions, then $\pi_{0}f(x)=x$; by the universality of the condition, then $f=\eta_{X}$.

    Now we show that (2) implies (3). Assume that $f:X\to V_{G}(X)$ is a coalgebra for this comonad. Let $U\subseteq X$ be a clopen subset. We want to show that
    \begin{equation*}
        {\uparrow}^{-1}[[U]]=\{ x : {\uparrow}x\in [U]\}=\Box U
    \end{equation*}
    is clopen. Consider the subset
    \begin{equation*}
        U^{*}\coloneqq \{z\in V_{G}(X) : z(0)\in U \text{ and } z(1)\in [U]\}.
    \end{equation*}
Note that this is a clopen subset of $V_{G}(X)$. Hence $f^{-1}[U^{*}]$ is a clopen subset of $X$.  Now note that $f(x)=(x,C_{0}^{x},C_{1}^{x},...)$. If $x\leq y$, then $f(x)\leq f(y)$, so since $C_{0}^{x}\leq C_{0}^{y}$, we have that $C_{0}^{y}\subseteq C_{0}^{x}$, and so, $y\in C_{0}^{x}$. In other words, $C_{0}^{x}={\uparrow}x$. 

Now then, note that
\begin{equation*}
    f^{-1}[U^{*}]=\Box U.
\end{equation*}
Indeed if $x\in \Box U$, then ${\uparrow}x\in [U]$, then by what we just discussed, $f(x)(0)\in U$ and $f(x)(1)\in [U]$, so by the above remark, $f(x)\in U^{*}$. Similarly, if $f(x)\in U^{*}$ it is clear to see that $x\in \Box U$. Hence because $f$ is continuous, we have that $\Box U$ is clopen, as desired.

Now we show that (3) implies (4). Note that this is the map sending a set $C$ to the set ${\uparrow}C$. Indeed if $U$ is clopen in $X$, then
\begin{equation*}
    {\uparrow}^{-1}[[U]]=\{C\in V(X) : {\uparrow}C\in [U]\}=\{ C\in V(X) : C\subseteq \Box U\}=[\Box U]
\end{equation*}
and we have that $\Box U$ is clopen, since it is exactly ${\uparrow}^{-1}[[U]]$ via the continuous map ${\uparrow}:X\to V(X)$. 

Finally, to show (4) implies (1), i.e., that $X$ is an Esakia space, assume that $U$ is clopen. We want to show that ${\downarrow}U$ is open. By pulling $\langle U\rangle$ along ${\uparrow}U$ we have that
\begin{equation*}
    {\uparrow}^{-1}[\langle U\rangle]=\{C : {\uparrow}C\cap U\neq \emptyset\}=\{C : C\cap {\downarrow}U\neq \emptyset\}=\langle {\downarrow}U\rangle.
\end{equation*}
Since the latter is open, by continuity of ${\uparrow}$, it can be written as
\begin{equation*}
    \bigcup_{i\in I}[V_{i}]\cap \langle W_{j_{0},i}\rangle \cap...\cap \langle W_{j_{k},i}\rangle
\end{equation*}
for clopens $V_{i}$ and $W_{j,i}$. Now if $x\in {\downarrow}U$, then $\{x\}\in \langle {\downarrow}U\rangle $, so $\{x\}\in [V_{i}]\cap \langle W_{j_{0},i}\rangle \cap...\cap \langle W_{j_{k},i}\rangle$ for some $j,i$; i.e., $V_{i}\cap W_{j_{0},i}\cap...\cap W_{j_{k},i}$ will be an open neighbourhood of $x$ entirely contained in ${\downarrow}U$. This shows that ${\downarrow}U$ is open, showing that $X$ is an Esakia space.\end{proof}

As a consequence we obtain that Esakia spaces can be seen, coalgebraically, as arising exactly from Priestley spaces by applying the $V_{G}$ construction. Moreover, we have the following:

\begin{theorem}
    The functor $V_{G}$ is comonadic.
\end{theorem}
\begin{proof}
    In light of Theorem \ref{comonadic adjunction}, all that is left to check is that p-morphisms correspond to coalgebra morphisms. If $X,Y$ are Esakia spaces, and $p:X\to Y$ is a p-morphism, we need to show that:
    \begin{equation*}
        V_{G}(p)\circ \eta_{X}=\eta_{Y}\circ p.
    \end{equation*}
    Note that these are two p-morphisms, and due to the universlity of the construction, two such p-morphisms will be equal if their projection on the first coordinate is the same; and indeed, we see that for $x\in X$:
    \begin{equation*}
        \pi_{0}\circ V_{G}(p)\circ \eta_{X}(x)=p(x) \text{ and } \pi_{0}\circ \eta_{Y}\circ p(x)=p(x).
    \end{equation*}
    Which shows that the desired diagram commutes.
    
    Conversely if $p:X\to Y$ is a coalgebra morphism, assume for $x\in X$ and $y\in Y$ that $p(x)\leq y$. Then $\eta_{Y}(p(x))\leq p(y)$, so $V_{G}(p)(\eta_{X}(x))\leq p(y)$. Hence by arguments similar to Lemma \ref{Lemma on lifting maps}, we can find $x\leq x'$ such that $p(x')=y$, using the fact that we can find such an $x'$ for finitely many coordinates. Hence $p$ is a p-morphism, as desired.   
\end{proof}

As discussed in  \cite{coalgebraicintuitionisticmodal}, this perspective can be taken further by looking at the ways in which coalgebraic representations for Priestley spaces can be lifted to coalgebraic representations for Esakia spaces, with applications in the study of modal Heyting algebras of various sorts\footnote{In abstract terms, one can use the adjunction $V_{G}\vdash I$ to lift any endofunctor $F:\mathbf{Pries}\to \mathbf{Pries}$ to an endofunctor $\overline{F}:\mathbf{Esa}\to \mathbf{Esa}$, and know that if $F$ is a comonad (e.g., represents a class of coalgebras dual to some appropriate class of modal distributive lattices), then $\overline{F}$ will again be a comonad -- indeed, $\overline{F}$ can be computed as the left Kan extension of $F$ along $V_{G}$. I am indebted to Alexander Kurz for this observation.}. This is currently an active direction, with many questions left open \cite{AlmeidaBezhanishviliDukic}.

\section{Free Heyting Extensions of Boolean Algebras and Inquisitive Logic}\label{Regular Heyting algebras and free Heyting extensions of Boolean algebras}

In this section we provide a non-trivial application of the tools developed in the previous sections. The problem we will be addressing concerns, categorically, another adjunction.

\subsection{Free Heyting Extensions of Boolean algebras}

Indeed, the categories $\mathbf{HA}$ and $\mathbf{BA}$, respectively of Heyting algebras with Heyting algebra homomorphisms and Boolean algebras with Boolean algebra homomorphisms, are related by a chain of adjunctions:
\begin{equation*}
\mathsf{F} \dashv \mathsf{Reg} \dashv I\dashv \mathsf{Center},
\end{equation*}
where
\begin{enumerate}
\item $I:\mathbf{BA}\to \mathbf{HA}$ is the inclusion;
\item $\mathsf{Center}:\mathbf{HA}\to \mathbf{BA}$ takes the \textit{center} of a Heyting algebra $H$, namely $\mathsf{Center}(H)=\{a\in H : a\vee \neg a=1\}$;
\item $\mathsf{Reg}:\mathbf{HA}\to \mathbf{BA}$ takes the \textit{regular elements} of a Heyting algebra $H$ namely $\mathsf{Reg}(H)=\{a\in H : \neg\neg a=a\}$.
\end{enumerate}
The functor $\mathsf{F}:\mathbf{BA}\to \mathbf{HA}$, in turn, is guaranteed to exist by the Adjoint Functor theorem, and to be the left adjoint of $\mathsf{Reg}$; a description using a word construction can be derived from the work of Moraschini \cite{Moraschini2018} (see also \cite[Example 4.18]{NicolauAlmeida2023}). It was likewise studied in \cite{Tur2008}, where it was shown to be fully faithful. We will show that such a functor can be described using the same ideas as above. 

We will prove some general facts about the construction before turning to some logical import of these facts. We will denote units and counits of the adjunctions at play with a specific superscript. Namely, when dealing with the natural transformations between functors between $\mathbf{Stone}$ and $\mathbf{Esa}$, we will denote by $\eta^{Sp}$ and $\varepsilon^{Sp}$ the proposed unit and counit; in turn, when dealing with the natural transformations between functors between $\mathbf{HA}$ and $\mathbf{BA}$, we will denote them as $\eta^{Alg}$ and $\varepsilon^{Alg}$.

\begin{lemma}\label{Continuity lemma}
    Let $X$ be an Esakia space, $Y$ a Stone space, and $f:\max X\to Y$ a continuous map. Then if $U\subseteq Y$ is clopen, $\neg\neg f^{-1}[U]$ is clopen in $X$.
\end{lemma}
\begin{proof}
    Indeed, note that $f^{-1}[U]$ is clopen in $\max X$, so there is some clopen $V$ such that $V\cap \max X=f^{-1}[U]$. Then note that $\neg\neg V=\neg\neg f^{-1}[U]$. Since $V$ is clopen, and the space $X$ is Esakia, this immediately gives us the result.
\end{proof}

Given $X$ a Stone space, let $e_{X}:\max V(X)\to X$ be the map which sends each singleton in $V(X)$ to its unique element; note that this is a natural isomorphism.

\begin{proposition}
    Assume that $X$ is an Esakia space, $Y$ is a Stone space and $f:\max X\to Y$ is a continuous map. Then there exists a unique order-preserving map $\tilde{f}:X\to V(Y)$, such that $e_{Y}\circ \tilde{f}\restriction_{\max}=f$ and $\tilde{f}$ is a p-morphism on maximal elements (i.e., if $\tilde{f}(x)\leq y$ and $y\in V(Y)$ is maximal, then there is some $w\geq x$ such that $\tilde{f}(w)=y$).
\end{proposition}
\begin{proof}
    For each $x\in X$, note that $f[{\uparrow}x\cap \max X]$ is a closed subset of $Y$. We claim that this provides the desired map. First note that it is continuous: for each $U\subseteq Y$ a clopen, we have that
    \begin{equation*}
        \tilde{f}^{-1}[[U]]=\{x : f[{\uparrow}x\cap \max X]\subseteq U\}=\{x : {\uparrow}x\cap \max X \subseteq f^{-1}[U]\}=\{x : x\in \neg\neg f^{-1}[U]\},
    \end{equation*}
    and by Lemma \ref{Continuity lemma}, the latter set is clopen; this is enough to show continuity.

    Note that then if $x\in \max X$, then $\tilde{f}(x)=f[{\uparrow}x\cap \max X]=\{f(x)\}$, and so by applying the counit we obtain the image of $f$, so the diagram commutes.

    It is easy to show that $\tilde{f}$ is order-preserving; we show it is unique. Indeed assume that $g:X\to V(Y)$ is another map forcing the diagram to commute. If $x\leq y$ and $y$ is maximal, then since $e_{Y}\circ g\restriction_{\max}(y)=f(y)$, and $g(x)\leq g(y)$, we have that $\{f(y)\}\subseteq g(x)$, i.e., $f(y)\in g(x)$. Conversely, if $y\in g(x)$, note that then $g(x)\leq \{y\}$, and this is maximal in $V(Y)$; hence there is some $w\geq x$ and $g(w)=\{y\}$. Let $w'\geq w$ be maximal. Then $g(w)\leq g(w')$, so $g(w')=\{y\}$ as well. Then $e_{Y}\circ g(w')=f(w')$, and so $w'\in {\uparrow}x\cap Max(X)$, and $f(w')=y$, so $y\in \tilde{f}(x)$. This shows that $\tilde{f}(x)=g(x)$ for each $x$.
\end{proof}

We now proceed by considering a variation of the $V_{g}(X)$ construction.

\begin{lemma}
    Let $X$ be a Stone space, and consider $V(X)$ the Vietoris space over $X$, understood as a Priestley space. then $V_{Max}(X)=\{C\in V_{r}(V(X)) : \forall D\in C, \forall x\in D, \{x\}\in C\}$ is a Priestley space.
\end{lemma}
\begin{proof}
    It suffices to show that it is a closed subset of $V_{r}(V(X))$, which is a Priestley space. Assume that $C\in V_{r}(V(X))$, and it is not in $V_{Max}(X)$; then let $D\in C$ and $x\in D$ be such that $\{x\}\notin C$. Then note that ${\uparrow}\{\{x\}\}\cap {\downarrow}C=\emptyset$, so in $V(X)$ we can find a clopen upset $V$ such that $\{\{x\}\}\subseteq V$ and $C\subseteq X-V$. Hence consider the neighbourhood:
    \begin{equation*}
        [X-V]\cap \langle {\downarrow}V\rangle
    \end{equation*}
    This is a clopen neighbourhood of $V_{Max}(X)$, since by Lemma \ref{Vietoris is Priestley}, $V(X)$ is Esakia. We have that $C$ belongs there, since it is entirely contained in $X-V$, and it contains $D$, and $D\leq \{x\}$ where $\{x\}\in V$. Moreover, if $E$ is any subset in this neighbourhood, $E\subseteq X-V$, and $E$ contains some subset $H$ where $H\leq H'$ and $H'\in V$. Let $H''\geq H'$ be maximal above $H'$ in $V(X)$. Then $H''$ is of the form $\{y\}$ for some $y\in X$, and $H''\in V$ because $V$ is an upset. Since $\{y\}\in V$, $\{y\}\notin E$, which shows the result.
\end{proof}

The restriction $r:V_{\max}(X)\to V(X)$ continues to be a surjective continuous function, since for each $C\in V(X)$ the principal upset ${\uparrow}C$ will be in $V_{\max}(X)$. It enjoys the following property:

\begin{lemma}
    Assume that $X$ is an Esakia space, $Y$ is a Stone space, and $f:\max X\to Y$ is a continuous map. Then there is a unique continuous and order-preserving $r$-open map $g_{f}:X\to V_{Max}(Y)$ with the property that $r\circ g_{f}=\tilde{f}$.
\end{lemma}
\begin{proof}
    We define $g_{f}:X\to V_{Max}(Y)$ by letting $g_{f}(x)=\tilde{f}[{\uparrow}x]$. This is certainly a rooted closed subset, and we show it has the characteristic property of belonging to the maximum-preserving Vietoris space. Indeed note that for each $x$:
    \begin{equation*}
        g_{f}(x)=\{f[{\uparrow}y\cap Max(X)] : x\leq y\}.
    \end{equation*}
    Now suppose that $D\in g_{f}(x)$. Then $D=f[{\uparrow}y\cap Max(X)]$ for $x\leq y$. Assume that $z\in D$, hence for some $w\geq y$, $f(w)=z$ where $w$ is maximal in $X$. Then $x\leq w$, so $\tilde{f}(w)\in g_{f}(x)$; but $\tilde{f}(w)=f[{\uparrow}w\cap Max(X)]=\{z\}$. This shows the property.

    The fact that the map is continuous and makes the diagram commute, and is unique, now follows from the same arguments from Lemma \ref{Duality Lemma for Key Property}. 
\end{proof}

Now assume that $f:\max X\to Y$ is a continuous function, where $X$ is an Esakia space and $Y$ is a Stone space. Define $M_{0}(Y)=V(Y)$, $M_{1}(Y)=V_{\max}(Y)$, $M_{n+1}(Y)=V_{r_{n}}(M_{n}(Y))$, where $r_{1}:V_{\max}(Y)\to V(X)$ is the root map, and $r_{n+1}:V_{n+1}(X)\to V_{n}(X)$ likewise. Let $M_{\infty}(Y)$ be the projective limit of this sequence.

\begin{lemma}
    Let $Y$ be a Stone space. Then $M_{\infty}(Y)$ is an Esakia space. Moreover, $\max(M_{\infty}(Y))\cong Y$.
\end{lemma}
\begin{proof}
    The fact that this is an Esakia space follows from noting that it is isomorphic to $V_{G}^{r_{1}}(M_{1}(Y))$. To see the fact about maximal elements, note that maximal elements of the projective limit will be sequences of maximal elements in $M_{n}(Y)$. Hence it suffices to show that for each of $M_{n}(Y)$, $\max(M_{n}(Y))\cong Y$. This is clear for $V(Y)$. For $V_{\max}(Y)$, note that if $C$ is maximal, since it is non-empty, there is some $D\in C$, and $x\in D$. By the condition, $\{x\}\in C$, so $C\leq \{\{x\}\}$. But then by maximality, $C=\{\{x\}\}$.

    Inductively, assume that $C\in M_{n+1}(Y)$, and $C$ is maximal. Look at $D\in M_{n}(Y)$ such that $D$ is maximal and $D\geq r_{n+1}(C)$; by induction hypothesis, $D=\{E\}$ for $E$ maximal in $M_{n-1}(Y)$. By assumption, then there is some $D'\in C$ such that $r_{n}(D)=r_{n}(D')=E$; but then $D'=\{E\}$. Then $C\leq {\uparrow}D'\cap C$, so by maximality of $C$ they are the same. Hence $r_{n+1}(C)=r_{n+1}({\uparrow}D\cap C)=D'=\{E\}$. This shows that $C$ is a singleton, as intended.
\end{proof}

Given a Stone space $Y$, let $\varepsilon_{Y}^{Sp}:\max(M_{\infty}(Y))\to Y$ be the natural isomorphism associating to each sequence of singletons $(x_{0},x_{1},...)$ where $x_{0}=\{x\}$ and $x$ is maximal in $Y$, the value $x$. It is easily checked that this is a natural isomorphism.

\begin{proposition}\label{Lifting of maps on maximal elements}
    Let $X$ be an Esakia space, $Y$ a Stone space, and $f:\max X\to Y$ a continuous map. Then there exists a unique order-preserving map $f^{\bullet}:X\to M_{\infty}(Y)$ such that $\varepsilon_{Y}\circ f^{\bullet}\restriction_{\max}=f$.
\end{proposition}
\begin{proof}
    Consider the lifting $f^{\bullet}$ defined by:
    \begin{enumerate}
        \item $f_{0}:X\to V(Y)$ is $\tilde{f}$ as defined above;
        \item $f_{1}:X\to V_{\max}(Y)$ is $g_{f}$ as defined above;
        \item $f_{n+1}:X\to M_{n+1}(Y)$ is given as follows from the $V_{G}^{g}$-construction, by successive direct images,
    \end{enumerate}
    and letting $f^{\bullet}(x)=(f_{0}(x),f_{1}(x),...)$. The fact that this is a continuous function follows from the universal property of projective limits, and it is certainly order-preserving. The argument for why it is a p-morphism is the same as in the case of the $V_{G}^{g}$-construction. Now finally, suppose that $x\in X$ is maximal; then since $f^{\bullet}$ is a p-morphism, $f^{\bullet}(x)$ is maximal as well, i.e., it is a sequence $(\{f(x)\},...)$. then $\varepsilon_{Y}\circ f^{\bullet}(x)=f(x)$ as desired.

    Now assume that $h:X\to M_{\infty}(Y)$ enjoys the same properties. Then first note that $h_{0}=\pi_{0}\circ h$ is a map $h_{0}:X\to V(Y)$, and $e_{Y}\circ h_{0}\restriction_{\max}=f$: indeed, if $x\in \max X$, then $h(x)$ is maximal, so $\pi_{0}\circ h(x)$ is a singleton $\{y\}$, and $y=f(x)$. Thus $h_{0}=\tilde{f}$. Similarly, we show that $h_{1}=g_{f}$, and successively, $h_{n+1}=f_{n}$, using the uniqueness properties. This shows that $h=f^{\bullet}$ as desired.
\end{proof}

The above allows us to define $M:\mathbf{Stone}\to \mathbf{Esa}$ a functor in a unique way: given a Stone space $Y$, $M(Y)\cong M_{\infty}(Y)$, and given a continuous function $f:Y\to Z$, we consider $f'=f\circ \pi_{0}:M_{\infty}(Y)\to Z$, and then lift uniquely as described by the above lemmas. We have thus shown:

\begin{theorem}\label{Adjunction of M and max}
    The functor $M:\mathbf{Stone}\to \mathbf{Esa}$ is right adjoint to $\max:\mathbf{Esa}\to \mathbf{Stone}$.
\end{theorem}

Algebraically, since $\max:\mathbf{Esa}\to \mathbf{Stone}$ is dual to the regularization functor \cite[A.2]{Esakiach2019HeyAlg}, we obtain:

\begin{corollary}
    The functor $M$ is dual to the functor $\mathsf{F}:\mathbf{BA}\to \mathbf{HA}$ which is left adjoint to $\mathsf{Reg}:\mathbf{HA}\to \mathbf{BA}$.
\end{corollary}

\subsection{Regular Heyting algebras and Inquisitve Logic}

More than simply giving an adjunction, the algebras $\mathsf{F}(A)$ actually appear naturally in the study of \textit{inquisitve logic}. Such a logic was extensively studied by Ciardelli \cite{Ciardelli2010}, having been introduced to model questions in classical logic, and having intimate connections to superintuitionistic logics. In \cite{Bezhanishvili2019_grilletti_holliday}, an algebraic approach to such logics was presented which gives them a semantics in terms of special Heyting algebras called \textit{regularly generated} Heyting algebras, or \textit{regular} for short:

\begin{definition}
    Let $H$ be a Heyting algebra. We say that $H$ is \textit{regular} if $H=\langle H_{\neg\neg}\rangle$, i.e., it is generated as a Heyting algebra by its regular elements. Let $\mathsf{RegHA}$ be the full subcategory of $\mathsf{HA}$ of regular Heyting algebras and Heyting homomorphisms.
\end{definition}

Indeed, in this paper the \textit{inquisitive extensions} arising there correspond to nothing more than the duals of the spaces $V(X)$ for $X$ a Stone space. As noted there, such algebras are always models of $\mathsf{Med}$, Medvedev's logic of finite problems. 

In \cite{bezhanishvili_grilletti_quadrellaro_2021} and \cite{Grilletti2023}, these Heyting algebras and their dual Esakia spaces were studied in connection with inquisitive logic. Here we will show that regular Heyting algebras actually arise in a categorically natural way, as those algebras for which the counit map is surjective. To prove these facts, we will need two preliminary lemmas:

\begin{lemma}\label{First step generation lemma}
    Let $X$ be a Stone space. Then $\mathsf{ClopUp}(V(X))$ is generated as a distributive lattice by its regular elements. 
\end{lemma}
\begin{proof}
    For each $U\subseteq X$ a clopen, notice that:
    \begin{equation*}
        \neg[U]=\{C\in V(X) : \forall D\subseteq C, D\nsubseteq U\}=\{C\in V(X) : C\subseteq X-U\}=[-U];
    \end{equation*}
    this follows because if $C\in \neg [U]$, then for each $x\in C$, $C\leq \{x\}$, so $\{x\}\nsubseteq U$ means that $x\in X-U$. This means that for each clopen $U$, $[U]$ is a regular clopen. The clopen upsets of $V(X)$ are unions of sets of this form, which gives us the result.
\end{proof}

\begin{lemma}\label{Second step generation lemma}
    If $X$ is a Stone space, the algebra $\mathsf{ClopUp}(V_{\max}(X))$ is generated from the elements of the form $[[U]]$ for $U\subseteq X$ a clopen, by unions and one layer of implications. 
\end{lemma}
\begin{proof}
    Notice that since $V(X)$ is generated by unions of sets of the form $[U]$ by Lemma \ref{First step generation lemma}, we have by the algebraic description of $\mathsf{ClopUp}(V_{r}(V(X)))$ that this algebra is unions of implications of unions of regular elements. Since $V_{\max}(V(X))$ is a subspace, it will be a distributive lattice quotient, and so the same will hold of this algebra.
\end{proof}

Hence we can prove:

\begin{theorem}\label{Freely generated Heyting algebra is regular}
    If $\alg{B}$ is a Boolean algebra, $\mathsf{F}(\alg{B})$ is a regular Heyting algebra.
\end{theorem}
\begin{proof}
    Dualizing the construction, we obtain a chain of injections:
    \begin{equation*}
        \alg{B}_{0}\rightarrow\alg{B}_{1}\rightarrow...\rightarrow\alg{B}_{n}\rightarrow ...\alg{B}_{\infty}
    \end{equation*}
    where $\mathsf{F}(\alg{B})=\alg{B}_{\infty}$ is a direct limit of this chain. Now let $a\in \alg{B}_{0}$ be a regular element. Then we claim that it remains regular in $\alg{B}_{\infty}$. To see this, the key observation is that in $V_{\max}(V(X))$, for $U\subseteq X$ a clopen:
    \begin{equation*}
        [V(X)-[U]]=[[X-U]],
    \end{equation*}
    which entails that the negation is preserved on the first step, and by the inductive construction, it will then be preserved at each subsequent step. If $C\in [V(X)-[U]]$, then $C\subseteq V(X)-[U]$, let $D\in C$. Want to show that $D\in [X-U]$. Indeed let $x\in D$. By our assumption, $\{x\}\in C$, so $\{x\}\notin [U]$, and so $x\in X-U$. Conversely if $C\in [[X-U]]$, then if $D\in C$, then $D\in [X-U]$, and so $D\notin [U]$ on pain of $D$ being empty.

    Hence we have that $\neg_{\alg{B}_{1}}a=\neg_{\alg{B}_{0}}a$. Consequently, if $a\in \alg{B}_{0}$ is regular, it remains regular in $\alg{B}_{1}$, and subsequently in all $\alg{B}_{n}$, meaning that $[a]\in \alg{B}_{\infty}$ is regular as well. By Lemmas \ref{First step generation lemma} and \ref{Second step generation lemma}, $\alg{B}_{0}$ as a subalgebra of $\alg{B}_{\infty}$ will be generated by regular elements as a distributive lattice, and $\alg{B}_{1}$ will be generated by $\alg{B}_{0}$ by implications of elements there, consequently, also by regular elements. Generally $\alg{B}_{n+1}$ will be generated by implications from elements of $\alg{B}_{n}$. Hence $\alg{B}_{\infty}$ will be generated as a Heyting algebra by its regular elements.
\end{proof}

\begin{definition}
    Let $B$ be a Boolean algebra. We call an algebra of the form $\mathsf{F}(B)$ a \textit{free regular Heyting algebra}.
\end{definition}

One can make such algebras appear naturally in the context of this construction\footnote{I am indebeted to Mat\'{i}as Menni for pointing out this fact, which allowed us to correct a previous, incorrect, formulation of the main result of this section.}:

\begin{proposition}
    The free regular Heyting algebras are precisely the coalgebras for the comonad induced by $\mathsf{F}\dashv\mathsf{Reg}$.
\end{proposition}
\begin{proof}
First note that it follows from our description of Theorem \ref{Adjunction of M and max} that the unit of the adjunction is an isomorphism; hence by \cite[Proposition 4.2.3]{Borceux1994}, it follows that $\mathsf{F}\dashv\mathsf{Reg}$ is idempotent. Hence, given any $\mathcal{B}$ a Boolean algebra, we have that:
\begin{equation*}
    \mathsf{F}(\mathcal{B})\cong \mathsf{F}(\mathsf{Reg}(\mathsf{F}(\mathcal{B}));
\end{equation*}
since algebras of the form $\mathsf{F}(\mathsf{Reg}(H))$ always admit a coalgebra structure, this ensures that each $\mathsf{F}(\mathcal{B})$ is a coalgebra for the comonad. Again by \cite[Proposition 4.2.3]{Borceux1994}, we have that all the coalgebras for this comonad must be isomorphisms, meaning that these are the only possible coalgebra structures.
\end{proof}

From this we will obtain a categorical description of regular Heyting algebras. We make use of the following observation:

\begin{lemma}\label{Surjective image of a regular Heyting algebra is regular}
    Let $H,H'$ be Heyting algebras, and suppose that $f:H\to H'$ is a surjective homomorphism. Suppose that $H$ is regularly generated. Then $H'$ is also regularly generated.
\end{lemma}
\begin{proof}
    For each $a\in H'$, there is some $b\in H$ such that $f(b)=a$. By regular generation, there is a Heyting algebra polynomial $\phi(x_{1},...,x_{n})$ and regular elements $b_{1},...,b_{n}$ such that $\phi(b_{1},...,b_{n})=b$. So $a=\phi(f(b_{1}),...,f(b_{n}))$. Note that for each regular element $b_{i}$, $\neg\neg f(b_{i})=f(\neg\neg b_{i})=f(b_{i})$, so $f(b_{i})$ are regular elements of $H'$. This shows that $H'$ is regularly generated.
\end{proof}

Now recall the structure of the counit of the adjunction. Dually, this is given as follows: given an Esakia space $X$, the unit $\eta_{X}^{Sp}:X\to M(\mathsf{max}(X))$, is obtained by lifting the identity on maximal elements as described by Proposition \ref{Lifting of maps on maximal elements}. Then note the following:

\begin{lemma}\label{Characterization of the counit and surjectivity}
    The counit $\varepsilon^{Alg}_{H}:\mathsf{F}(\mathsf{Reg}(H))\to H$ is surjective if and only $H$ is dual to a regularly generated Heyting algebra.
\end{lemma}
\begin{proof}
    By Theorem \ref{Freely generated Heyting algebra is regular} and  Lemma \ref{Surjective image of a regular Heyting algebra is regular}, if $\varepsilon^{Alg}_{H}$ is surjective, then $H$ is regular. Conversely, assume that $H$ is regular, and let $X$ be its dual Esakia space. Let $\eta_{X}^{Sp}(x)=\eta_{X}^{Sp}(y)$ for $x,y\in X$. Since $x\neq y$, without loss of generality, assume that $x\nleq y$, and let $U$ be a clopen upset separating them. By assumption, $U$ is generated from the regular elements, say $U=\phi(V_{0},...,V_{n})$ where $\phi$ is some Heyting polynomial and $V_{i}$ are regular clopen upsets. We now show by induction that for any $w,v\in X$, if $\eta_{X}^{Sp}(w)=\eta_{X}^{Sp}(v)$, then for any $\phi(\overline{V})$ a polynomial over regular clopen upsets, $w\in \phi(\overline{V})$ if and only if $v\in \phi(\overline{V})$:
    \begin{enumerate}
        \item First assume that $U=V_{i}$.  Since $\eta_{X}^{Sp}(w)=\eta_{X}^{Sp}(v)$, we have that $\max({\uparrow}w)=\max({\uparrow}v)$; since $w\in U$, then $\max({\uparrow}w)\subseteq U$, so $\max({\uparrow}v)\subseteq U$, and so, since $U$ is regular, $v\in U$.
        \item Inductively, assume that for each $\phi(x_{1},...,x_{n})$ of implication-depth at most $n$, $w\in \phi(V_{0},...,V_{n})$ if and only if $v\in \phi(V_{0},...,V_{n})$ for $V_{i}$ regular clopen upset. It is clear that this condition persists when taking intersections and unions. So assume that $\psi=\phi_{0}\rightarrow \phi_{1}$. Let $\overline{V}=(V_{0},...,V_{n})$ be a tuple of regular clopen upsets. Assume that $w\notin \psi(\overline{V})$. then $w\leq w'$ and $w'\in \phi_{0}(\overline{V})$ and $w'\notin \phi_{0}(\overline{V})$. Then $\eta_{X}^{Sp}(w)\leq \eta_{X}^{Sp}(w')$, so since $\eta_{X}^{Sp}(w)=\eta_{X}^{Sp}(v)$, and the fact that this is a p-morphism, there is some $v\leq v'$ such that $\eta_{X}^{Sp}(w')=\eta_{X}^{Sp}(v')$. By induction hypothesis, $w'\in \phi_{0}(\overline{V})-\phi_{1}(\overline{V})$ entails that $v'\in \phi_{0}(\overline{V})-\phi_{1}(\overline{V})$, and so $v\notin \psi(\overline{V})$.
    \end{enumerate}
    We thus conclude that assuming that $\eta_{X}^{Sp}(x)=\eta_{X}^{Sp}(y)$, and that $x\neq y$, we obtain a contradiction; hence $x=y$, meaning that $\eta_{X}^{Sp}$ is injective. By duality, this means that $\varepsilon_{H}^{Alg}$ is surjective, as desired.
\end{proof}

From this we obtain the following description:

\begin{theorem}
    The full subcategory $\mathsf{RegHA}$ coincides precisely with the full subcategory $$\mathsf{Epi}(\mathsf{F})=\{H \in \mathsf{HA} : \varepsilon_{H}^{Alg} \text{ is a (regular) epimorphism}\}.$$
\end{theorem}
\begin{proof}
    This follows immediately from Lemma \ref{Characterization of the counit and surjectivity}. 
\end{proof}

The above answers, in a categorical fashion, a question left open in \cite{Grilletti2023} as to the structure of Esakia duals of regularly generated Heyting algebras.

\section{Subvarieties of Heyting algebras}\label{Subvarieties of Heyting algebras}

Unlike products, coproducts of algebras can vary drastically when one moves from a variety $\mathbf{K}$ to a subvariety $\mathbf{K}'$. As such, in this section, we briefly discuss with a few examples how the above construction can be adapted to handle some well-known subvarieties of $\mathbf{HA}$. We will begin by outlining how this works for Boolean algebras, which will allow us to recover the well-known construction of the \textit{free Boolean extension} of a distributive lattice \cite{gehrkevangoolbooktopologicalduality}. We will then provide analogous free constructions of $\mathsf{KC}$-algebras and $\mathsf{LC}$-algebras\footnote{Since completing this paper we have learned of independent work by Carai \cite{caraifreegodelalgebras} outlining the construction of Free G\"{o}del algebras. For a comparison with that approach one can check Section 6 of that paper, which explains the relationship between Carai's construction and ours.}; we leave a detailed analysis of the full scope of the method for future work.

$\mathsf{KC}$-algebras are Heyting algebras axiomatised with an additional axiom
\begin{equation*}
    \neg p \vee \neg\neg p;
\end{equation*}
they are also sometimes called WLEM-algebras (For Weak Excluded Middle) or DeMorgan Algebras. They appear prominently in settings like Topos Theory (see e.g. \cite{caramellotopostheoretic}). In turn, $\mathsf{LC}$-algebras are Heyting algebras axiomatised by
\begin{equation*}
    p\rightarrow q\vee q\rightarrow p;
\end{equation*}
they are sometimes called ``G\"{o}del algebras" or ``G\"{o}del-Dummett algebras" as well as ``prelinear Heyting algebras" (the latter is due to the fact that the subdirectly irreducible elements in such a variety are chains). They are of key interest in the study of many-valued logic and substructural logic \cite[Chapter 3]{Hjek1998}.

\subsection{Boolean algebras}

Our general method in handling subvarieties will be to take the Ghilardi complex construction, and show that by quotienting by additional equations at  appropriate levels, one can obtain sub-complexes which limit to the free algebra in the desired variety. To illustrate this, we consider the simple case of Boolean algebras: let $\alg{D}$ be a distributive lattice. Consider $\mathbf{F}_{DLat}(\{(a,b) : a,b\in D\})/\Theta\cup \{p\vee \neg p=1\}$, the free algebra quotiented as before, with the addition of the equation $p\vee \neg p$, with which one forces Booleannness. Dually, given $X$ a Priestley space, this takes a given subspace of the space $V_{1}(X)$ (recall, by definition \ref{Ghilardi complex definition} that this is simply the set of closed and rooted subsets, $V_{r}(X)$):

\begin{lemma}\label{Booleanisation from step-by-step lemma}
    Let $X$ be a Priestley space. Given a point $C\in V_{1}(X)$, we have that $C\in [\neg U]\cup [U]$ for each clopen upset $U\subseteq X$ if and only if $C$ is a singleton. 
\end{lemma}
\begin{proof}
    If $C$ is a singleton this is clear. Now assume that $C$ has two points $x,y$. Then either $x\nleq y$ or $y\nleq x$; without loss of generality assume the former. Then $x\in U$ and $y\notin U$ for some clopen upset. So $C\notin [U]$ and $C\notin [\neg U]$, as desired.
\end{proof}

Let $B_{1}(X)\subseteq V_{1}(X)$ be the set of singletons of $X$. Using Lemma \ref{Booleanisation from step-by-step lemma} one can prove the following:

\begin{theorem}
    For each Priestley space $X$, $B_{1}(X)\cong X$ as Stone spaces, and $B_{1}(X)$ is isomorphic to the Priestley space $(X,=)$.
\end{theorem}
\begin{proof}
    The isomorphism is clear; the freeness follows by the same arguments as in Lemma \ref{Lifting Lemma on the Esakia spaces}.
\end{proof}

We thus get that $B_{1}(X)$, the one-step extension from distributive lattices to Boolean algebras, actually coincides with the well-known free Boolean extension, the left adjoint to the inclusion of Boolean algebras in distributive lattices.

It is also worthy to note that Boolean algebras have a peculiar feature as far as our construction $V_{G}$ goes: if $X$ is a Priestley space with a discrete order (definitionally equivalent to a Stone space), then $V_{G}(X)\cong X$, and in fact $V_{1}(X)\cong X$ already (since rooted subsets of singletons have to simply be singletons as well). Stone spaces also exhibit some further interesting properties in this respect:

\begin{theorem}\label{Stabilisation Theorem for the construction}
    Let $X$ be a Priestley space; then the following are equivalent:
    \begin{enumerate}
        \item $X$ is a Stone space;
        \item $V_{1}(X)\cong X$ where the root map witnesses the isomorphism;
        \item For some $n$, $V_{n}(X)\cong V_{n+1}(X)$ where the root map witnesses the isomorphism.
        \item $V_{G}(X)\cong X$ through the projection $\pi_{0}$, i.e., the counit of the adjunction is an isomorphism for $X$.
    \end{enumerate}
\end{theorem}
\begin{proof}
By the previous remark (1) implies (2), and similar considerations mean that (2) implies (3). We show that (3) implies (1).

    Assume that $X$ is not a Stone space, i.e., there is a pair of points $x,y$ such that $x\leq y$. Note that then $x<y$ is a subposet of $X$. Recall from Example \ref{step by step construction of Rieger-Nishimura} that from such a poset one constructs the whole Rieger-Nishimura ladder, and there is throughout each layer of the construction a point $z_{n}$ such that $r$ does not provide an isomorphism. Since all sets involved are finite, this means that $V_{n}(X)$ will always likewise contain such a point $z_{n}$ where $r$ cannot provide an isomorphism.

    Next note that certainly (1) implies (4). Now suppose that $X$ is not a Stone space. Then as noted $V_{1}(X)\not\cong X$, so there are at least two elements $C_{0},C_{1}\in V_{1}(X)$ which have the same root. But then we can construct extensions of such elements to elements of $V_{G}(X)$ whilst sharing the same root, i.e., the projection from $V_{G}(X)$ is not injective.
\end{proof}

The above theorem does not rule out that there can be other Priestley spaces $X$ such that $V_{G}(X)\cong X$, but it implies that such isomorphisms will not be natural with respect to the adjunction. Such a question is left as an interesting open problem.

\subsection{KC-algebras}

To generalise our construction to $\mathsf{KC}$ we will work with the axiomatization of this logic given by adding to $\mathsf{IPC}$ the axiom $\neg p\vee \neg\neg p$. Since  this equation has rank $2$, our restriction will need to happen at the second stage, which is where one can appropriately define the law.

\begin{definition}
    Let $g:X\to Y$ be a continuous and order-preserving map. Let $C\in V_{r_{g}}(V_{g}(X))$ be an element of the second stage of the $g$-Ghilardi complex over $X$. We say that $C$ is \textit{well-directed} if whenever $D,D'\in C$, we have that ${\uparrow}D\cap {\downarrow}D'\neq \emptyset$.
    We write $V_{g,K}(V_{g}(X))$ for the subset of $V_{r_{g}}(V_{g}(X))$ consisting of the well-directed elements.
\end{definition}

The motivation for this definition is the following, where we recall that given $g:X\to Y$ a Priestley morphism, $V_{2}(X)=V_{r_{g}}(V_{g}(X))$:

\begin{proposition}\label{Duality for free KC-algebras}
    Let $g:X\to Y$ be a continuous and order-preserving function. Given $C\in V_{2}(X)$ we have that $C\in [[X-U]]\cup [V_{1}(X)-[X-U]]$ for each clopen upset $U\subseteq X$ if and only if $C$ is well-directed.
\end{proposition}
\begin{proof}
    Assume that  $C\notin [V_{1}(X)-[X-U]]\cup [[X-U]]$. In other words, there is some $D\in C$ such that $D\notin V_{1}(X)-[X-U]$, and some $D'\notin [X-U]$. The former means that $D\subseteq X-U$, and the later means that $D'\cap U\neq \emptyset$. Let $x\in D'\cap U$. Then $D''={\uparrow}x\cap D$ is such that $D'\leq D''$, and so because $C$ is $r_{g}$-open, there is some $E\in C$ such that $D'\leq E$ and $r_{g}(D'')=r_{g}(E)$, i.e., $x$ is the root of $E$. Since $U$ is a clopen upset, $E\subseteq U$, and $D\subseteq -U$, so certainly
    \begin{equation*}
        E\cap U=\emptyset.
    \end{equation*}
    This shows that $C$ is not well-directed.

    Conversely, assume that $C$ is not well-directed. Let $D,D'\in C$ be such that ${\uparrow}D\cap {\downarrow}D'=\emptyset$. Using the properties of Priestley spaces, namely Strong Zero-Dimensionality, we can then find a clopen upset $U$ such that $D\subseteq U$ and $D'\subseteq X-U$; but this means precisely that $C\notin [[X-U]]\cup [V_{1}(X)-[X-U]]$, by the arguments above.
\end{proof}

Proposition \ref{Duality for free KC-algebras} shows that the construction $V_{g,K}$ indeed yields a Priestley space, and so, that we can proceed inductively. In order to conclude our results, we will need to show, however, that this construction is free. Recall that we say that a Priestley space $X$ is \textit{directed} if whenever $x,y,z\in X$, and $x\leq y,z$ there is some $y,z\leq w$. Then we can show:

\begin{lemma}\label{Inductive lemma for KC algebras}
    Let $Z$ be a Priestley space which is directed, and let $h:Z\to X$ be a $g$-open, continuous and order-preserving map. Let $h^{*}:Z\to V_{g}(X)$ be the canonical map yielded by Lemma \ref{Key Lemma in Universal Construction}.  Then the unique $r_{r}$-open, continuous and order-preserving map $h'$ making the obvious triangle commute factors through $V_{g,K}(V_{g}(X))$.
\end{lemma}
\begin{proof}
    Similar to Lemma \ref{Key Lemma in Universal Construction} we have that the definition of $h^{**}$ is entirely forced. All that we need to do is verify that given $x\in Z$, $h^{**}$ is well-directed.

    For that purpose, note that
    \begin{equation*}
        h^{**}(x)=h^{*}[{\uparrow}x]=\{h^{*}(y) : x\leq y\}
    \end{equation*}
    Now assume that $x\leq y$ and $x\leq y'$. Because $Z$ is directed, there is some $z$ such that $y\leq z$ and $y'\leq z$. Then note that
    \begin{equation*}
        h^{*}(z)\subseteq h^{*}(y)\cap h^{*}(y'),
    \end{equation*}
    since if $z\leq m$, then $y\leq m$ and $y'\leq m$, so $h(m)\in h^{*}(y)\cap h^{*}(y')$. This shows that $h^{**}(x)$ is well-directed, as desired.
\end{proof}

Denote by $V_{GK}(X)$ the projective limit of the Ghilardi complex obtained by iterating the construction $V_{g,K}$, as before. Then using these constructions we obtain:

\begin{theorem}\label{Key Theorem on free Godel algebras}
Let $\alg{A}$ be a Heyting algebra, let $\alg{D}$ be a distributive lattice, and let $p:\alg{A}\to \alg{D}$ be a distributive lattice homomorphism which preserves the relative pseudocomplements from $\alg{A}$ (such that $\alg{D}$ contains them). Let $\mathcal{K}$ be the dual of $V_{GK}(X)$. Then $\mathcal{K}$ is the $\mathsf{KC}$-algebra freely generated by $\alg{D}$ which preserves the relative pseudocomplements coming from $p$.
\end{theorem}
\begin{proof}
    By our earlier remarks, we have that $\mathcal{K}$ is the directed union of Heyting algebras, and must satisfy the Weak Excluded Middle. To show freeness, suppose that $\alg{M}$ is a $\mathsf{KC}$ algebra and $k:\alg{D}\to \alg{M}$ is a map preserving the relative pseudocomplements coming from $\alg{A}$. Then using Lemma \ref{Inductive lemma for KC algebras} systematically, we see that $\alg{M}$ will factor through $K_{C}^{n}(H(\alg{D}))$ for every $n$, and hence, a fortriori, through the limit. This shows the result.
\end{proof}

Note that the above construction is not guaranteed to provide an \textit{extension} of the algebra $\alg{D}$ -- it might be that some elements are identified by the first quotient. However, it is not difficult to show that if $\alg{D}$ is dual to a directed Priestley space, then it embeds into $\alg{K}(\alg{D})$ (by showing, for instance, that the composition with the root maps to $\alg{D}$ is always surjective). We could also run similar constructions as those done for Heyting algebras in the previous sections, with a completely parallel theory.

\subsection{LC-Algebras}

We finally consider the case of $\mathsf{LC}$ algebras. As is well-known, such a variety is locally finite, which means that we ought to expect our step-by-step construction to in fact stabilise at a given finite step when the original algebra is finite. As we will see, much more is true. Let us say that a Priestley space $X$ is \textit{prelinear} if for each $x\in X$, ${\uparrow}x$ is a linear order.

\begin{definition}
    Let $g:X\to Y$ be a monotone and continuous function. We say that $C\in V_{g}(X)$ is a \textit{linearised} closed, rooted and $g$-open subset, if $C$ is a chain. Denote by $V_{g,L}(X)$ the set of linearised closed, rooted and $g$-open subsets.
\end{definition}

As before, we can prove: 

\begin{proposition}
    Let $g:X\to Y$ be a continuous and order-preserving function. Then for each $C\in V_{1}(X)$, we have that $C\in [-U\cup V]\cup [-V\cup U]$ if and only if $C$ is linearised.
\end{proposition}
\begin{proof}
    Assume that $U,V$ are clopen subsets of $X$, and $C\notin  [-U\cup V]\cup [-V\cup U]$. This amounts to saying that
    \begin{equation*}
        C\nsubseteq -U\cup V \text{ and } C\nsubseteq -V\cup U.
    \end{equation*}
    In other words, $x\in C$ and $x\in U-V$ and $y\in C$ and $y\in V-U$. Since these sets are clopen subsets, this means that $x\nleq y$ and $y\nleq x$.

    Conversely, assume that $C$ is not linearised. Since $x\nleq y$, there is some clopen upset $U$ such that $x\in U$ and $y\notin U$, and since $y\nleq x$, there is some clopen upset $V$ such that $y\in V$ and $x\notin V$. These two subsets then witness that $C\notin  [-U\cup V]\cup [-V\cup U]$.
\end{proof}

We now prove, again, an analogue to Lemma \ref{Key Lemma in Universal Construction}:

\begin{lemma}
Let $Z$ be a prelinear Priestley space, and let $h:Z\to X$ be a $g$-open, continuous and order-preserving map. Then the unique $r_{g}$-open, continuous and order-preserving map $h'$, making the obvious triangle commute, factors through $V_{g,L}(X)$.  
\end{lemma}
\begin{proof}
    Simply note that since $Z$ is prelinear, we have that $h[{\uparrow}x]$ is linear, given that ${\uparrow}x$ is as well.
\end{proof}

Denote by $V_{G}^{L}(X)$ the projective limit of the sequence 
\begin{equation*}
    (X,V_{1}^{L}(X),V_{2}^{L}(X),...)
\end{equation*}
where $V_{i}^{L}(X)$ is defined as in Definition \ref{Ghilardi complex definition} except that from the second step onwards we use $V_{g,L}(-)$ rather than $V_{g}(-)$. Then in the same way as before, we can prove:

\begin{theorem}
Let $\alg{A}$ be a Heyting algebra, let $\alg{D}$ be a distributive lattice, and let $p:\alg{A}\to \alg{D}$ be a distributive lattice homomorphism which preserves the relative pseudocomplements from $\alg{A}$ (such that $D$ contains them). Let $\mathcal{K}$ be the dual of $V_{g,L}^{\infty}(X)$. Then $\mathcal{L}$ is the $\mathsf{LC}$-algebra freely generated by $D$ which preserves the relative pseudocomplements coming from $p$.
\end{theorem}

Consider a finite poset $P$. Because of local finiteness, we know that $V_{G}^{L}(P)\cong V_{n}^{L}(P)$ must hold for some finite $n$ (otherwise the construction would add infinitely many new points). By Theorem \ref{Stabilisation Theorem for the construction}, $V^{L}_{1}(X)\cong X$ can never hold for such a poset unless it is discrete, since $V_{1}^{L}(-)=V_{1}(-)$. But it is reasonable to ask whether $V_{n}^{L}(X)\cong V_{n+1}^{L}(X)$ holds for any poset. Notably, we will see that in fact the construction always stabilises in two steps, even for infinite spaces $X$.

\begin{lemma}\label{Lemma on comparable prelinear elements}
Let $C,D\in V_{L}^{2}(X)$. Suppose that $r(C)=r(D)$ and that $C$ and $D$ are comparable. Then $C=D$.
\end{lemma}
\begin{proof}
Assume that $x\in V_{L}^{2}(X)$, and let $R$ be its root. Suppose that $A\in x$ is arbitrary. Then we claim that
\begin{equation*}
    A={\uparrow}r(A)\cap R.
\end{equation*}
Indeed, if $b\in A$, then since $R\leq A$, $b\in R$, and by definition $b\in {\uparrow}r(A)$. Conversely, if $b\in R\cap {\uparrow}r(A)$, then $R\leq {\uparrow}b\cap R$ so there is some $B\in x$ such that $r(B)=b$. Note that since $x$ is a chain, we have that $A\leq B$ or $B\leq A$; in the first case $b\in A$, in the latter, since $r(A)\leq b$, we have $r(A)=b$, which again means that $b\in A$. In either case we obtain equality.

Now assume that $C,D\in x$ and $r(C)=r(D)$. Then by what we have just showed:
\begin{equation*}
    C={\uparrow}r(C)\cap R={\uparrow}r(D)\cap R=D.
\end{equation*}
\end{proof}

\begin{proposition}\label{After two times and prelinearity we have stability}
Let $X$ be a Priestley space. Assume that $V_{2}^{L}(X)$ is prelinear. Then $V_{3}^{L}(X)\cong V_{2}^{L}(X)$.
\end{proposition}
\begin{proof}
We prove that the only rooted linearised and $r$-open subsets of $V_{2}^{L}(X)$ are of the form ${\uparrow}x$ for $x\in V_{2}^{L}(X)$. Let $C\subseteq {\uparrow}x$ be a linearised rooted and $r$-open subset, such that $r(C)=x$. If $C\neq {\uparrow}x$ there is some $y$ such that $x\leq y$ and $y\notin C$. Because $C$ is $r$-open, there is some $y'\in C$ such that $r(y)=r(y')$. Now note that $y'$ and $y$ are comparable (they are both successors of $x$), so by Lemma \ref{Lemma on comparable prelinear elements}, we have that $y=y'$, a clear contradiction. So we must have no such $y$ must exist, and hence, that no $C$ in these conditions can exist. We then conclude that every rooted linearised subset is of the form ${\uparrow}x$ for some $x$, which yields an easy isomorphism.
\end{proof}

\begin{proposition}\label{Applying the construction twice is enough}
Let $X$ be a Priestley space. Then $V_{2}^{L}(X)$ is prelinear.
\end{proposition}
\begin{proof}
Assume that $x\leq y$ and $x\leq z$, and look at $r(y)$ and $r(z)$. Since these are in $x$, which is a chain, assume that $r(y)\leq r(z)$. Then we claim that $y\leq z$. Indeed, assume that $A\in z$. Then $r(z)\leq A$, so $r(y)\leq A$. Because $y$ is $r$-open, let $A'$ be such that $A'\in y$, and $r(A)=r(A')$. Then $A'\in x$, and since $A\in z$, $A\in x$ as well; so by Lemma \ref{Lemma on comparable prelinear elements}, we have that $A=A'$. This shows that $A\in y$, which shows the result.
\end{proof}

Putting these Propositions together, we obtain the following result:

\begin{theorem}
Let $X$ be a Priestley space. Then $V_{2}^{L}(X)\cong V_{G}^{L}(X)$.
\end{theorem}
\begin{proof}
By Proposition \ref{Applying the construction twice is enough} we know that $V_{L}^{2}(X)$ will be prelinear; and by Proposition \ref{After two times and prelinearity we have stability} we have that iterating the construction past the second stage does not generate anything new. Hence for each $n>2$, $V_{n}^{L}(X)\cong V_{2}^{L}(X)$, meaning that the projective limit will be likewise isomorphic.
\end{proof}

The above theorem, coupled with some of the theory we have developed for Heyting algebras, and which could easily be generalised for G\"{o}del algebras, can have some computational uses: given two G\"{o}del algebras, to compute their coproduct, we can look at their dual posets, form the product, and apply $V_{L}^{2}$, to obtain the resulting product. This provides an alternative perspective from  the recursive construction discussed in \cite{DAntona2006,aguzzoligodelcoproducts}, and can be compared with the recent work of Carai in \cite{caraifreegodelalgebras}. 

In logical terms, this has some interesting consequences. Recall from \cite{Shehtman2016}\footnote{There this property is studied under the name of ``finite formula depth property"; the present name was introduced since it makes obvious the connection between such a property and local tabularity, as well as allowing simpler aliases, e.g., $\mathsf{LC}$ being a $2$-uniform logic.} that a superintuitionistic logic $L$ is said to be \textit{n-uniform} if for each formula $\phi$, $\phi$ is equivalent over $L$ to a formula with implication rank $n$. We say that $L$ is \textit{uniformly locally tabular} if there is some $n$ such that $L$ is $n$-uniform. It was shown by Shehtman that uniform local tabularity entails local tabularity, and it is open whether the converse holds. Our result above implies the following:

\begin{theorem}
    The logic $\mathsf{LC}$ is $2$-uniform.
\end{theorem}

This result invites the question of how one should convert an arbitrary formula into one of the normal forms of $\mathsf{LC}$. This sort of computational procedure requires a tighter understanding of the algebraic structure of the step by step construction, and we leave it for future work.

\section{Categories of Posets with P-morphisms}\label{The Category of Posets with P-morphisms}

We now turn our attention to how the questions from the previous sections can be understood in a ``discrete" setting -- i.e., by considering posets, rather than Esakia spaces. This can be motivated from a categorical point of view, where it arises as a natural question concerning the existence of adjunctions in $\mathsf{Ind}$ and $\mathsf{Pro}$-completions of categories. Another motivation for this comes from the logical setting, where such discrete duals correspond to the ``Kripke semantics" of modal and superintuitionistic logics. As we will see, a construction analogous to the $V_{G}$ functor can be provided, relating the category of \textit{Image-Finite Posets} with p-morphisms, and the category of posets with monotone maps.

\subsection{Pro- and Ind-Completions of categories of finite objects}

To explain the questions at hand, we recall the picture, outlined in \cite[Chapter VI]{johnstone1982stone}, relating categories of finite algebras and several of their completions; we refer the reader to \cite{Massas2023} where these and other dualities are outlined. To start with, we have $\mathbf{FinBA}$, the category of finite Boolean algebras, and its dual, $\mathbf{FinSet}$. For such categories, there are two natural ``completions":
\begin{itemize}
    \item The \textit{Ind-completion}, $\mathsf{Ind}(-)$, which freely adds all directed colimits;
    \item The \textit{Pro-completion}, $\mathsf{Pro}(-)$, which freely adds all codirected limits.
\end{itemize}
For categories of algebras, the Pro-completion can equivalently be seen as consisting of the profinite algebras of the given type, whilst the Ind-completion can be seen as obtaining all the \textit{locally finite} algebras of that same type -- which in the case of Boolean algebras or Distributive lattices, will be all algebras. As such, we have that:
\begin{equation*}
    \mathsf{Ind}(\mathbf{FinBA})=\mathbf{BA} \text{ and } \mathsf{Pro}(\mathbf{FinBA})=\mathbf{CABA}
\end{equation*}
where $\mathbf{BA}$ is the category of all Boolean algebras with Boolean algebra homomorphisms, and $\mathbf{CABA}$ is the category of complete and atomic Boolean algebras with complete homomorphisms. Correspondingly, we have that
\begin{equation*}
    \mathsf{Ind}(\mathbf{FinSet})=\mathbf{Set} \text{ and } \mathsf{Pro}(\mathbf{FinSet})=\mathbf{Stone}.
\end{equation*}
The facts that $\mathbf{Set}^{op}\cong \mathbf{CABA}$ (Tarski duality) and $\mathbf{Stone}^{op}\cong \mathbf{BA}$ (Stone duality) then amount to the basic relationship between Ind and Pro-completions, namely, for any category $\mathbf{C}$:
\begin{equation*}
    (\mathsf{Ind}(\mathbf{C}^{op}))^{op}\cong \mathbf{Pro}(\mathbf{C}) \text{ and } (\mathsf{Pro}(\mathbf{C}^{op}))^{op}\cong \mathsf{Ind}(\mathbf{C}).
\end{equation*}

A similar picture can be drawn for distributive lattices: there the relevant category of finite objects is the category $\mathbf{FinPos}$ of finite posets with monotone maps, which is dual to $\mathbf{FinDL}$, finite distributive lattices with lattice homomorphisms, and we have that
\begin{align*}
    \mathsf{Ind}(\mathbf{FinDL})=\mathbf{DL} &\text{ and } \mathsf{Pro}(\mathbf{FinDL})=\mathbf{CCJDL}\\
    \mathsf{Ind}(\mathbf{FinPos})=\mathbf{Pos} &\text{ and } \mathsf{Pro}(\mathbf{FinPos})=\mathbf{Pries},
\end{align*}
where $\mathbf{DL}$ is the category of all distributive lattices with lattice homomorphisms; $\mathbf{CCJDL}$ is the category of complete and completely join-prime generated distributive lattices with complete lattice homomorphisms; $\mathbf{Pos}$ is the category of all posets with monotone maps.

If we are interested in Heyting algebras, it is not hard to see that the category of finite Heyting algebras with Heyting algebra homomorphisms is dual to the category $\mathbf{FinPos}_{p}$ of finite posets with p-morphisms. The Ind-completion of $\mathbf{FinHA}$ is not necessarily $\mathbf{HA}$, since not every Heyting algebra is a directed colimit of finite Heyting algebras; however, one may use this intuition to ``guess" a discrete category which allows the generalization of the above results. Indeed, the natural place to look is the Ind-completion of $\mathbf{FinPos}_{p}$. Such objects are those which are directed colimits of finite posets via p-morphisms -- but this means that the finite posets must embed into the top of the object, i.e., such objects must be image-finite posets, in other words
\begin{equation*}
    \mathsf{Ind}(\mathbf{FinPos}_{p})\cong \mathbf{ImFinPos}
\end{equation*}
where $\mathbf{ImFinPos}$ is the category of image-finite posets with p-morphisms; this is, indeed, the main result from \cite{Bezhanishvili2008}, where it is shown that these are dual to the category of profinite Heyting algebras, as desired. Indeed, a similar picture to this one has been discussed for finite modal algebras, and their duals, finite Kripke frames, in the recent work \cite{profinitenessmonadicityuniversalmodels}, where the above dualities are briefly sketched.

As a result of this relationship, it is natural to look for an adjunction holding between the inclusions of $\mathbf{ProHA}$ into $\mathbf{CCJDL}$, or equivalently, for a right adjoint to the inclusion of $\mathbf{ImFinPos}$ into $\mathbf{Pos}$. This is what we outline in the next section.

\subsection{An Adjunction between the category of Image-Finite posets and Posets with Monotone Maps}

Our constructions in this section mirror very much the key constructions from Section \ref{Free Constructions of Heyting Algebras over Distributive Lattices}. Indeed, given a poset $g:P\to Q$, denote by
\begin{equation*}
    P_{g}(X)=\{C\subseteq X : C \text{ is rooted, $g$-open, and finite }\}.
\end{equation*}
As usual when $g$ is the terminal map to the one point poset, we write $P_{r}$. Also note that regardless of the structure of $X$, $P_{g}(X)$ will always be an image-finite poset, though in general the root maps will not admit any sections.

\begin{definition}
    Let $P,Q$ be posets with $g:P\to Q$ a monotone map. The \textit{g-Powerset complex} over $P$, $(P^{g}_{\bullet}(X),\leq_{\bullet})$ is a sequence
    \begin{equation*}
        (P_{0}(X),P_{1}(X),...)
    \end{equation*}
    connected by monotone maps $r_{i}:P_{i+1}(X)\to P_{i}(X)$ such that:
    \begin{enumerate}
        \item $P_{0}(X)=X$;
        \item $r_{0}=g$;
        \item For $i\geq 0$, $P_{i+1}(X)\coloneqq P_{r_{i}}(P_{i}(X))$;
        \item $r_{i+1}\coloneqq r_{r_{i}}:P_{i+1}(X)\to P_{i}(X)$ is the root map.
    \end{enumerate}
    We write $P_{G}^{g}(X)$ for the \textit{image-finite part} of the projective limit of the above sequence in the category of posets with monotone maps. In other words, we take
    \begin{equation*}
        \{x\in \varprojlim P_{n}(X) : {\uparrow}x \text{ is finite }\}.
    \end{equation*}
\end{definition}

We now prove the following, which is analogous to the property from Lemma \ref{Duality Lemma for Key Property}:

\begin{lemma}\label{Duality Lemma for Key Property on Posets}
    Let $g:X\to Y$ be a monotone map between posets; $h:Z\to X$ be a monotone and $g$-open map, where $Z$ is image-finite. Then there exists a unique monotone and $r_{g}$-open map such that the triangle in Figure \ref{fig:commutingtriangleofposets} commutes. 

\begin{figure}[h]
    \centering
\begin{tikzcd}
Z \arrow[rd, "h"'] \arrow[rr, "h'"] &   & P_{g}(X) \arrow[ld, "r"] \\
& X &                         
\end{tikzcd}    \caption{Commuting Triangle of Posets}
    \label{fig:commutingtriangleofposets}
\end{figure}
\end{lemma}
\begin{proof}
    The arguments will all be the same except we now need to show that given $a\in Z$, $h'(a)$ is a finite subset; but since we assume that ${\uparrow}a$ is finite, this immediately follows.
\end{proof}

Using this we can show that the poset $P_{G}(X)$ will be free in the desired way. The key fact that is missing concerns the fact that, after applying Lemma \ref{Duality Lemma for Key Property on Posets} infinitely often, one needs, as in Proposition \ref{Lifting Lemma on the Esakia spaces}, to show that the lifting is indeed a p-morphism. There we used compactness, which certainly does not hold for arbitrary posets, and this is the place where image-finiteness becomes crucial.

\begin{proposition}
    Given a poset $X$, $P_{G}^{g}(X)$ has the following universal property: whenever $p:Z\to X$ is a monotone map from an image-finite poset $Z$, there is a unique extension of $p$ to a p-morphism $\overline{p}:Z\to P_{G}(X)$.
\end{proposition}
\begin{proof}
    Using Lemma \ref{Duality Lemma for Key Property on Posets} repeatedly we construct a map $f:Z\to \varprojlim P_{n}(X)$. We will show that such a map is in fact a p-morphism; if we do that, then since $Z$ is image-finite, it will factor through the image-finite part of the projective limit, and hence will give us the desired map.

    So assume that $x\in Z$, and suppose that
    \begin{equation*}
        f(x)\leq y
    \end{equation*}
    Note that by construction, $y=(a_{0},a_{1},...)$ for some elements, sending the root map appropriately. Now let $n$ be arbitrary. Then consider $a_{n+1}$, which by definition is a subset of $p_{n+1}[{\uparrow}x]=\{p_{n}(y) : x\leq y\}$. That means that there is some $y'$ such that $x\leq y'$ and $f(y')$ agrees with $y$ up to the $n$-the position (which follows from the commutation of the triangles in the above Lemma). Since $x\in Z$ has only finitely many successors, this means that there must be a successor $x\leq a$, such that $f(a)$ agrees with $y$ on arbitrarily many positions, i.e., $f(a)=y$. This shows the result as desired.
\end{proof}

We can then show:

\begin{proposition}
    The map $P_{G}:\mathbf{Pos}\to \mathbf{ImFinPos}$ is a functor; indeed it is left adjoint to the inclusion of $I:\mathbf{ImFinPos}\to \mathbf{Pos}$.
\end{proposition}
\begin{proof}
    One verifies, analogously to Lemma \ref{Lemma on lifting maps} that the construction $P_{G}^{g}$ is functorial, and the rest follows just like in Theorem \ref{Adjunction Theorem}.
\end{proof}

A few facts are worthy of note here: following the main ideas of \cite{profinitenessmonadicityuniversalmodels}, the above describes an adjunction which splits the adjunction between profinite Heyting algebras and the category of sets. Moreover, as noted in \cite{profinitenessmonadicityuniversalmodels}, the construction $P_{G}$ given above can be thought of as a generalization of the $n$-universal model -- indeed, if one starts with the dual poset of the free distributive lattice on $n$-generators, $P_{G}$ will produce precisely this same model. This illustrates an interesting connection between these two well-known constructions of the free algebra -- the Ghilardi/Urquhart step-by-step construction, and the Bellissima/Grigolia/Shehtman universal model -- showing that they are, in some sense, dual to each other.

It is also notable that the above methods cannot be expected to be extended to yield an adjunction between $\mathbf{Pos}_{p}$ -- the category of posets with p-morphisms -- and $\mathbf{Pos}$ as above. Indeed, as noted in \cite{openmapsdoesnothaveproducts}, $\mathbf{Pos}_{p}$ does not have binary products; on the other hand, if $I:\mathbf{Pos}_{p}\to \mathbf{Pos}$ had a right adjoint, then all limits in $\mathbf{Pos}$ would have to be preserved.

We conclude by noting that the two constructions provided -- $V_{G}$ and $P_{G}$ -- are indeed intimately related for any finite poset $P$. For this we will need a few more technical developments:

\begin{definition}
    Let $(V_{n}(X))_{n\in \omega}$ be a Ghilardi complex. Let $x\in V_{n}(X)$ be any point. We say that $x$ is prestable if $r_{n+1}^{-1}[x]$ is a singleton; we say that it is stable if whenever $x\leq x'$, then $x'$ is prestable\footnote{These concepts seem to have been first considered by Dito Pataraia in unpublished work. I credit Pataraia for the main ideas of these results, though the proofs are slightly different.}.
\end{definition}

\begin{lemma}\label{Stable points remain stable}
Assume that $x\in X_{n+1}$. If $r_{n+1}(x)$ is stable, then $x$ is stable as well.
\end{lemma}
\begin{proof}
Assume that $C,D$ are such that $r_{n+2}(C)=r_{n+2}(D)=x$. Suppose that $y\in C$. Then $x\leq y$. Since $D$ is $r_{n+1}$-open, there is some $y'\in D$ such that $r_{n+1}(y)=r_{n+1}(y')$. Since $r_{n+1}(x)$ is stable, then $y=y'$, so $y\in D$. This shows that $C=D$, which, since ${\uparrow}x$ roots at $x$, proves that $x$ is prestable. The same argument shows that all successors of $x$ are prestable, obtaining the result.
\end{proof}

\begin{lemma}\label{Pataraia stability lemma}[Pataraia's Stability Lemma]
Let $x\in V_{n}(X)$ be arbitrary. Then $x^{\bullet}=\{{\uparrow}({\uparrow}y) : x\leq y\}\in X_{n+2}$ is stable, and $r_{n+1}(r_{n+2}(x))=x$.
\end{lemma}
\begin{proof}
Note that $x^{\bullet}$ is $r_{n+1}$-open. Now assume that $r_{n+3}(C)=r_{n+3}(D)=x^{\bullet}$. Let $Y\in C$; then $x^{\bullet}\leq Y$. This means that $r_{n+2}(Y)={\uparrow}({\uparrow}y)$ for $x\leq y$. Now suppose that $y\leq z$. Then note that ${\uparrow}({\uparrow}z)\geq r_{n+2}(Y)$; since the latter is $r_{n+1}$-open, there is some $K\in Y$ such that $r_{n+1}(K)={\uparrow}z$. Note that $Y'=Y\cap {\uparrow}K$ is then such that $x^{\bullet}\leq Y'$, so by the same argument, $K={\uparrow}({\uparrow}m)$ for some $m$; but then given its root, $K={\uparrow}({\uparrow}z)$. Hence $Y=y^{\bullet}$. 

Note that by $r_{n+2}$-openness of $D$, there is some $Y'\in D$ such that $r_{n+2}(Y')=r_{n+2}(Y)$. The exact same argument then shows that $Y'=y^{\bullet}$. So $Y\in D$, which shows that $C=D$. This shows that $x^{\bullet}$ is prestable. Now note that then, by the argument we have just shown, any successor of $x^{\bullet}$ is of the form $y^{\bullet}$, so the former is stable.
\end{proof}

\begin{corollary}\label{Infinite stabilisation}
    If $n\in \omega$ is arbitrary, and $x\in V_{n}(X)$, there is some point $x'\in V_{G}(X)$ such that ${\uparrow}x'$ is finite, and $x'(n)=x$.
\end{corollary}
\begin{proof}
    Using Lemma \ref{Pataraia stability lemma} and Lemma \ref{Stable points remain stable}, we can construct an extension of $x$ to a point of $V_{G}(X)$ which only has finitely many extensions, using the fact that the root map is an isomorphism on the restriction from $x^{\bullet}$ to $x$, and onwards along the construction.
\end{proof}

\begin{theorem}\label{Profiniteness of the image finite}
    Let $P$ be a finite poset. Then the image finite part of $V_{G}(P)$ is dense in the space.
\end{theorem}
\begin{proof}
For this we need to show that for each $U\subseteq V_{n}(P)$ a (clopen) set, $\pi_{n}^{-1}[U]$ contains a point $x$ which only sees finitely many points. Indeed, suppose that $x(n)\in U$. Using Corollary \ref{Infinite stabilisation}, let $\overline{x}$ be a finite point such that $\overline{x}(n)=x(n)$. Then $\overline{x}\in \pi_{n}^{-1}[U]$. This shows that the image-finite part is dense, as desired.
\end{proof}

Recall that given a poset $(P,\leq)$, an ordered-topological space $(X,\leq,\tau)$ and an order-preserving map $p:P\to X$, we say that $p$ is an \textit{order-compactification} if $p$ is an order-homeomorphism onto the image, and $p[P]$ is dense in $X$ (see \cite{nachbin1965topology}). Then we have:

\begin{corollary}
    For each finite poset $P$, $V_{G}(P)$ is an order-compactificationof $P_{G}(P)$.
\end{corollary}
\begin{proof}
    By Theorem \ref{Profiniteness of the image finite} this follows immediately.
\end{proof}

\section{Conclusions and Further Research}\label{Conclusion}

In this paper we have explicitly described the adjunction between Heyting algebras and Distributive lattices, and have extracted some consequences for the theory of Heyting algebras from this -- namely, facts about the structure of free Heyting algebras, properties of the category of Heyting algebras, and direct proofs of amalgamation of Heyting algebras deriving from the corresponding properties for distributive lattices. We have also looked at how this construction fares in different settings: when considering specific subvarieties of Heyting algebras (Boolean algebras, $\mathsf{KC}$-algebras, $\mathsf{LC}$-algebras), where it is shown that the construction can be reasonably adapted, and sometimes corresponds to natural modifications of the Ghilardi functor; and also we considered the relationship between the category of image-finite posets with p-morphisms and the category of posets with monotone maps, showing that reasonable modifications also yield an adjunction here.

As noted throughout the text, the material presented here does not exhaust the range of applications of these constructions and the ideas wherein: in ongoing work \cite{AlmeidaBezhanishviliDukic}, the results of the present paper have applications in the study of coalgebraic logic; in future work we also expect to be able to use these ideas to investigate the structure of normal forms in superintuitionistic logics, as well as issues such as the existence and complexity of interpolants. Moreover, the results of this paper leave open several questions, both of a technical and of a conceptual nature.

One fact which we have not worked with concerns the structure of the Priestley spaces obtained by applying $V_{g}(-)$; when $X$ i a finite poset, this yields a finite poset, which has a natural structure as a (bi-)Esakia space; when $X$ is linear, then $V_{r}(X)$ is homeomorphic to $V(X)$, and hence by the results of Section \ref{Regular Heyting algebras and free Heyting extensions of Boolean algebras}, is an Esakia space. But it is not clear that even starting with a bi-Esakia space $X$, $V_{g}(X)$ will be bi-Esakia.

Following in the study of subvarieties of Heyting algebras, a systematic study connecting correspondence of axioms with Kripke semantics, and the appropriate modifications made to the $V_{G}$ functor, seems to be in order. This does not seem immediately straightforward, as our example with $\mathsf{KC}$ algebras illustrates, but it might be possible to obtain more general results, at least in the case of finitely axiomatisable such logics.

It would also be interesting to study similar constructions to the ones presented here for other signatures; a natural extension would be to study bi-Heyting algebras and bi-Esakia spaces. We expect that this should provide a fair challenge, since rather than simply adding relative pseudocomplements, such a construction would need to also add relative supplements, and even in the finite case, no such description appears to be available in the literature.

\section{Acknowledgements}

I would like to thank, in no special order, Silvio Ghilardi, Guram Bezhanishvili, Nick Bezhanishvili, Marco Abbadini, Luca Carai, Qian Chen and Lingyuan Ye for several inspiring conversations and suggestions on parts of this work. Special thanks for Frederik Lauridsen and Matías Menni for pointing out imprecisions and incorrect statements in previous versions of this work. A deep thanks to Mamuka Jibladze for many enthusiastic discussions on the main construction of this paper, and for his key idea in proving the main lemma concerning the construction $V_{r}$.

\printbibliography[
    heading=bibintoc,
    title={Bibliography}
]

\end{document}